\documentclass[12pt,reqno]{amsart}
\usepackage{amsmath , amssymb }
\usepackage{graphics}

\usepackage[mathscr]{eucal}

\newdimen\unit\newdimen\psep\newcount\nd\newcount\ndx\newbox\dotb\newbox\ptbox
\newdimen\dx\newdimen\dy\newdimen\dxx\newdimen\dyy\newdimen\hgt
\newdimen\xoff\newdimen\yoff
\newcommand\clap[1]{\hbox to 0pt{\hss{#1}\hss}}
\newcommand\vdisk[1]{{\font\dotf=cmr10 scaled #1\dotf.}}
\newcommand\varline[2]{\setbox\dotb\hbox{\vdisk{#1}}\xoff=-.5\wd\dotb
\wd\dotb=0pt\yoff=-.5\ht\dotb\psep=#2\ht\dotb}
\newcommand\varpt[1]{\setbox\ptbox\clap{\vdisk{#1}}\setbox\ptbox
\hbox{\raise-.5\ht\ptbox\box\ptbox}}
\newcommand\cpt{\copy\ptbox}
\newcommand\point[3]{\rlap{\kern#1\unit\raise#2\unit\hbox{#3}}}
\newcommand\setnd[4]{\dx=#3\unit\advance\dx-#1\unit\divide\dx by\psep
\dy=#4\unit\advance\dy-#2\unit\divide\dy by\psep \multiply\dx
by\dx\multiply\dy by\dy\advance\dx\dy\nd=1\advance\dx-1sp
\loop\ifnum\dx>0\advance\dx-\nd sp\advance\nd1\advance\dx-\nd
sp\repeat}
\newcommand\dl[4]{{\setnd{#1}{#2}{#3}{#4}\dline{#1}{#2}{#3}{#4}\nd}}
\newcommand\dline[5]{{\nd=#5\hgt=#2\unit\dx=#3\unit\advance\dx-#1\unit
\divide\dx by\nd\dy=#4\unit\advance\dy-#2\unit\divide\dy by\nd
\advance\hgt\yoff\rlap{\kern#1\unit\kern\xoff\loop\ifnum\nd>1\advance\nd-1
\advance\hgt\dy\kern\dx\raise\hgt\copy\dotb\repeat}}}

\newcommand\qellip[4]{{\setnd{0}{0}{#3}{#4}\dx=\unit\dy=0pt\raise\yoff\rlap{%
\kern#1\unit\kern\xoff\raise#2\unit\hbox{\loop\ifnum\dx>0\rlap{\kern#3\dx
\raise#4\dy\copy\dotb}\hgt=\dx\divide\hgt
by\nd\advance\dy\hgt\hgt=\dy \divide\hgt
by\nd\advance\dx-\hgt\repeat\rlap{\raise#4\dy\copy\dotb}}}}}
\newcommand\bez[6]{{\setnd{#1}{#2}{#3}{#4}\ndx=\nd\setnd{#3}{#4}{#5}{#6}
\ifnum\ndx>\nd\nd=\ndx\fi\dx=#3\unit\advance\dx-#1\unit\dy=#4\unit
\advance\dy-#2\unit\dxx=#5\unit\advance\dxx-#1\unit\dyy=#6\unit\advance
\dyy-#2\unit\advance\dxx-2\dx\advance\dyy-2\dy\divide\dxx
by\nd\divide\dyy
by\nd\advance\dx.25\dxx\advance\dy.25\dyy\divide\dx
by\nd\divide\dy by\nd \multiply\nd
by2\dx=100\dx\dy=100\dy\dxx=100\dxx\dyy=100\dyy\divide\dxx by\nd
\divide\dyy
by\nd\hgt=#2\unit\raise\yoff\rlap{\kern#1\unit\kern\xoff
\raise\hgt\copy\dotb\loop\ifnum\nd>0\advance\nd-1\advance\hgt0.01\dy
\kern0.01\dx\raise\hgt\copy\dotb\advance\dx\dxx\advance\dy\dyy\repeat}}}
\newcommand\ptu[3]{\point{#1}{#2}{\cpt\raise1ex\clap{$\scriptstyle{#3}$}}}
\newcommand\ptd[3]{\point{#1}{#2}{\cpt\raise-1.8ex\clap{$\scriptstyle{#3}$}}}
\newcommand\ptr[3]{\point{#1}{#2}{\cpt\raise-.4ex\rlap{$\ \scriptstyle{#3}$}}}
\newcommand\ptl[3]{\point{#1}{#2}{\cpt\raise-.4ex\llap{$\scriptstyle{#3}\ $}}}
\newcommand\ptlu[3]{\point{#1}{#2}{\raise.8ex\clap{$\scriptstyle{#3}$}}}
\newcommand\ptld[3]{\point{#1}{#2}{\raise-1.6ex\clap{$\scriptstyle{#3}$}}}
\newcommand\ptlr[3]{\point{#1}{#2}{\raise-.4ex\rlap{$\,\scriptstyle{#3}$}}}
\newcommand\ptll[3]{\point{#1}{#2}{\raise-.4ex\llap{$\scriptstyle{#3}\,$}}}
\newcommand\pt[2]{\point{#1}{#2}{\cpt}}
\newcommand\thkline{\varline{1600}{.3}}
\newcommand\medline{\varline{800}{.5}}
\newcommand\thnline{\varline{400}{.6}}

\varpt{2500}\thnline\unit=2em

\newtheorem{thm}{Theorem}
\newtheorem*{thmA}{Theorem A}                                                            
\newtheorem{conj}{Conjecture}

\newtheorem{prob}{Problem}
\newtheorem{lemma}[thm]{Lemma}
\newtheorem{cor}[thm]{Corollary}

\newtheorem{obs}[thm]{Observation}
\theoremstyle{definition}\newtheorem{rmk}{Remark}
\theoremstyle{definition}
\newcommand{\ds}{\displaystyle}
\newcommand{\ul}{\underline}

\def\A{\mathcal{A}}
\def\B{\mathcal{B}}
\def\C{\mathcal{C}}
\def\D{\mathcal{D}}

\def\F{\mathcal{F}}
\def\G{\mathcal{G}}
\def\HH{\mathcal{H}}
\def\J{\mathcal{J}}

\def\P{\mathcal{P}}
\def\Q{\mathcal{Q}}
\def\R{\mathcal{R}}
\def\S{\mathcal{S}}

\def\AA{\mathbb{A}}

\def\N{\mathbb{N}}
\def\QQ{\mathbb{Q}}
\def\RR{\mathbb{R}}
\def\Z{\mathbb{Z}}

\def\le{\leqslant}
\def\ge{\geqslant}

\def\eps{\varepsilon}

\begin{document}
\title[Ordered graphs]{Hereditary properties of ordered graphs}

\author{J\'ozsef Balogh}
\address{Department of Mathematics\\ University of Illinois\\ 1409 W. Green Street\\ Urbana, IL 61801} \email{jobal@math.uiuc.edu}

\author{B\'ela Bollob\'as}
\address{Department of Mathematical Sciences\\ The University of Memphis\\ Memphis, TN 38152\\ and\\ Trinity College\\ Cambridge CB2 1TQ\\ England} \email{bollobas@msci.memphis.edu}

\author{Robert Morris}
\address{Department of Mathematical Sciences\\ The University of Memphis\\ Memphis, TN 38152} \email{rdmorrs1@memphis.edu}\thanks{The first author was supported during this research by OTKA grant T049398 and NSF grant DMS-0302804, the second by NSF grant ITR 0225610, and the third by a Van Vleet Memorial Doctoral Fellowship.}

\begin{abstract}
An ordered graph is a graph together with a linear order on its vertices. A hereditary property of ordered graphs is a collection of ordered graphs closed under taking order-preserving isomorphisms of the vertex set, and order-preserving induced subgraphs. If $\P$ is a hereditary property of ordered graphs, then $\P_n$ denotes the collection $\{G \in \P : V(G) = [n]\}$, and the function $n \mapsto |\P_n|$ is called the speed of $\P$.

The possible speeds of a hereditary property of labelled graphs have been extensively studied (see~\cite{BBW1} and \cite{ICM} for example), and more recently hereditary properties of other combinatorial structures, such as oriented graphs (\cite{AlekS}, \cite{BBMtou}), posets (\cite{BBMpos}, \cite{BGP}), words (\cite{BBwor}, \cite{QZ}) and permutations (\cite{KK}, \cite{MT}), have also attracted attention. Properties of ordered graphs generalize properties of both labelled graphs and permutations.

In this paper we determine the possible speeds of a hereditary property of ordered graphs, up to the speed $2^{n-1}$. In particular, we prove that there exists a jump from polynomial speed to speed $F_n$, the Fibonacci numbers, and that there exists an infinite sequence of subsequent jumps, from $p(n)F_{n,k}$ to $F_{n,k+1}$ (where $p(n)$ is a polynomial and $F_{n,k}$ are the generalized Fibonacci numbers) converging to $2^{n-1}$. Our results generalize a theorem of Kaiser and Klazar~\cite{KK}, who proved that the same jumps occur for hereditary properties of permutations.
\end{abstract}

\maketitle

\section{Introduction}\label{orderintro}

In this paper we shall determine the possible speeds of a hereditary property of ordered graphs, up to the speed $2^{n-1}$. In particular, we shall prove that there is a jump from polynomial speed to speed $F_n$, the Fibonacci numbers, and that there exists an infinite sequence of subsequent jumps converging to $2^{n-1}$. Our results generalize a theorem of Kaiser and Klazar~\cite{KK}, who proved that the same jumps occur for hereditary properties of permutations. We begin by making the definitions necessary in order to state our main result.

An \emph{ordered graph} is a graph together with a linear order on its vertices. As a convention, we shall assume that if $G$ is an ordered graph of order $n$, then $V(G) = [n]$, where $i < j$ as vertices of $G$ if $i < j$ in $\N$. A collection of ordered graphs is called a \emph{property} if it is closed under order-preserving isomorphisms of the vertex set. Given ordered graphs $G$ and $H$, we say that $G$ is an \emph{induced ordered subgraph} of $H$ (and write $G \le H$) if there exists an injective, order-preserving map $\phi : V(G) \to V(H)$ such that $ij \in E(G)$ if and only if $\phi(i)\phi(j) \in E(H)$. A property of ordered graphs $\P$ is called \emph{hereditary} if it is closed under taking induced ordered subgraphs. In this paper `subgraph' will always mean `induced ordered subgraph', unless otherwise stated. Given a property of ordered graphs, $\P$, we write $\P_n$ for the collection of ordered graphs in $\P$ with vertex set $[n]$. The \emph{speed} of $\P$ is simply the function $n \mapsto |\P_n|$. Analogous definitions can be made for other combinatorial structures (e.g., graphs, posets, permutations).

We are interested in the (surprising) phenomenon, observed for hereditary properties of various types of structure (see for example~\cite{BBW1}, \cite{BGP}, \cite{MT}) that the speeds of such a property are far from arbitrary. More precisely, there often exists a family $\F$ of functions $f : \N \to \N$ and another function $F : \N \to \N$ with $F(n)$ {\em much} larger than $f(n)$ for every $f \in \F$, such that if for each $f \in \F$ the speed is infinitely often larger than $f(n)$, then it is also larger than $F(n)$ for every $n \in \N$. Putting it concisely: the speed {\em jumps} from $\F$ to $F$.

We can now state our main result. Let $F_{n,k}$ denote the $n^{th}$ generalized Fibonacci number of order $k$, defined by $F_{n,k} = 0$ if $n < 0$, $F_{0,k} = 1$ and $F_{n,k} = F_{n-1,k} + F_{n-2,k} + \ldots + F_{n-k,k}$ for every $n \ge 1$. We shall sometimes write $F_n$ for $F_{n,2}$.

\begin{thm}\label{order}
If $\P$ is a hereditary property of ordered graphs, then one of the following assertions holds.

\begin{enumerate}
\item[(a)] $|\P_n|$ is bounded, and there exist $M,N \in \N$ such that $|\P_n| = M$ for every $n \ge N$.\\[-1ex]
\item[(b)] $|\P_n|$ is a polynomial. There exist $k \in \N$ and integers $a_0, \ldots, a_k$ such that, $|\P_n| = \sum_{i = 0}^k a_i {n \choose i}$ for all sufficiently large $n$, and $|\P_n| \ge n$ for every $n \in \N$.\\[-1ex]
\item[(c)] $F_{n,k} \le |\P_n| \le p(n)F_{n,k}$ for every $n \in \N$, for some $2 \le k \in \N$ and some polynomial $p$, so in particular $|\P_n|$ is exponential.\\[-1ex]
\item[(d)] $|\P_n| \ge 2^{n-1}$ for every $n \in \N$.\\[-1ex]
\end{enumerate}
\end{thm}

We remark that our theorem is inspired in part by the work of Kaiser and Klazar~\cite{KK}, who proved an identical theorem for hereditary properties of permutations. In fact Theorem~\ref{order} is a generalization of their result, since every hereditary property of permutations $\Q$ may be thought of as a hereditary property of ordered graphs $\P(\Q)$ in the following way. For each $n \in \N$ and $\pi \in \Q_n$, let $G(\pi)$ be the ordered graph with vertex set $[n]$, and with edge set $\{ij :$ the order of the vertices $i$ and $j$ is reversed by $\pi\}$. Let $\P(\Q) = \{G(\pi) : \pi \in \Q\}$. It is easy to see that $\P(\Q)$ is hereditary, and that $|\P(\Q)_n| = |\Q_n|$. In fact, writing $\Pi$ for the collection of all permutations, one can give an even simpler description of the property $\P(\Pi)$. Let $H_1$ denote the ordered graph on vertex set $[3]$ and with edge set $\{12,23\}$, and $H_2$ that with edge set $\{13\}$. Then $\P(\Pi) = \{G : H_1 \not\le G$ and $H_2 \not\le G\}$, where $H_i \not\le G$ means that $H_i$ is not an induced ordered subgraph of $G$. Hence the theorem of Kaiser and Klazar is Theorem~\ref{order} in the case that $H_1,H_2 \notin \P$.

There are other interesting special cases of Theorem~\ref{order}. For example, let $\G$ be a hereditary property of oriented graphs, and let $\G^n_{mon}$ denote the collection of pairs $(G,\phi)$, where $G \in \G$, $|G| = n$ and $\phi: [n] \leftrightarrow V(G)$ is a monotone labelling of the vertices of $G$, i.e., if $x \to y$ in $G$ then $\phi(x) < \phi(y)$. Since each monotone labelling may be thought of as an ordering, there is a hereditary property of ordered graphs $\P$ such that $|\G^n_{mon}| = |\P_n|$ for every $n \in \N$. Hence the speed $n \mapsto |\G^n_{mon}|$ satisfies the conclusion of Theorem~\ref{order}.

Similarly, let $\R$ be a hereditary property of posets, and let $\P^n_{lin}$ denote the collection of pairs $(P,\psi)$, where $P \in \R$, $|P| = n$, and $\psi$ is a linear extension of $P$, i.e., a monotone labelling of the elements of $P$. Each pair $(P,\psi)$ may be thought of as a transitive monotone-labelled oriented graph, so the speed $n \mapsto |\R^n_{lin}|$ also satisfies the conclusion of Theorem~\ref{order}. We suspect that our list of interesting special cases is not exhaustive.

A property of \emph{graphs} is a collection of (unlabelled) graphs closed under isomorphism, and a property of graphs is hereditary if it is closed under taking (non-ordered) induced subgraphs. Given a property of graphs, $\P$, we write $\P^n$ for the collection of labelled graphs of order $n$ in $\P$ (i.e., the collection of non-isomorphic pairs $(G,\phi)$, where $G \in \P$, $|G| = n$ and $\phi: [n] \leftrightarrow V(G)$ is a labelling of the vertices of $G$), and call the function $n \mapsto |\P^n|$ the labelled speed of $\P$. The labelled speed of a property of graphs  was introduced by Erd\H{o}s~\cite{E} in 1964, and subsequently studied by Erd\H{o}s, Kleitman and Rothschild~\cite{EKR}, Erd\H{o}s, Frankl and R\"odl~\cite{EFR}, Kolaitis, Pr\"omel and Rothschild~\cite{KPR}, and Pr\"omel and Steger~\cite{PS3}, \cite{PS4}, \cite{PS5}, amongst others, though always in the special case where only a single graph is forbidden. The study of the possible speeds of a \emph{general} hereditary property of labelled graphs was initiated by Scheinerman and Zito~\cite{SZ} in 1994. They were the first to study speeds below $n^n$, and proved that for such properties the speeds all lie in a few fairly narrow ranges. A little later, considerably stronger results were proved by Balogh, Bollob\'as and Weinreich~\cite{BBW1}, \cite{BBW4}. In the range $2^{cn^2}$, the main results were proved by Alekseev~\cite{Alekseev}, by Bollob\'as and Thomason~\cite{BTbox}, \cite{BT}, and by Pr\"omel and Steger~\cite{PS5}. Putting all these results together, one obtains the following powerful theorem. Here $B_n$ denotes the $n^{th}$ Bell number, the number of partitions of $[n]$.

\begin{thmA}\label{thmA}
Let $\P$ be a hereditary property of graphs. Then one of the following holds.
\begin{enumerate}
\item[(a)] $|\P^n| = \ds\sum_{i=1}^k p_i(n) i^n$ for every $n \ge N$, for some $N, k \in \N$, and some collection $p_1(n), \ldots, p_k(n)$ of polynomials.\\[-1ex]

\item[(b)] $|\P^n| = n^{(1 - 1/k + o(1))n}$ as $n \to \infty$, for some $2 \le k \in \N$. \\[-1ex]

\item[(c)] $n^{(1+o(1))n} = B_n \le |\P^n| \le 2^{o(n^2)}$ as $n \to \infty$.\\[-1ex]

\item[(d)] $|\P^n| = 2^{(1 - 1/k + o(1)){n \choose 2}}$, as $n \to \infty$, for some $2 \le k \in \N$.\\[-1ex]

\item[(e)] $|\P^n| = 2^{n \choose 2}$ for every $n \in \N$.\\[-1ex]
\end{enumerate}
\end{thmA}

Given a property of graphs $\G$, one can define a property of ordered graphs $\P(\G)$ by taking every possible ordering on the vertex set of each graph in $\G$. Note that now $|\P(\G)_n| = |\G^n|$, so the speed of a property of ordered graphs is also a generalization of the labelled speed of a property of graphs. Using this idea, we can easily deduce the following theorem from Theorem A.

\begin{thm}\label{highord}
Let $\P$ be a hereditary property of ordered graphs. Either $|\P_n| = 2^{o(n^2)}$ as $n \to \infty$, or $|\P_n| = 2^{(1 - 1/k + o(1)){n \choose 2}}$ for some $2 \le k \in \N$.
\end{thm}

\begin{proof}
Given a hereditary property of ordered graphs, $\P$, we can naturally associate a property of graphs $\G$ with $\P$, by identifying isomorphic graphs with different linear orders. Since $\P$ is hereditary, so is $\G$. Also, since each linear order may be thought of as a labelling, and there are at most $n!$ different labellings of a graph in $\G_n$, we have
\begin{equation} |\G_n| \le |\P_n| \le |\G^n| \le n! \cdot {|\G_n|} \label{highineq} \end{equation}
for every $n \in \N$. So if $|\P_n| \ge 2^{cn^2}$ for some $c > 0$ and for infinitely many $n$, then also $|\G^n| \ge 2^{cn^2}$ for these $n$. Hence, by Theorem A, $|\G^n| = 2^{(1 - 1/k + o(1)){n \choose 2}}$ for some $2 \le k \in \N$.

Now, by \eqref{highineq}, we also have $|\G_n| = 2^{(1 - 1/k + o(1)){n \choose 2}}$, since $n! = 2^{o(n^2)}$, and so $|\P_n| = 2^{(1 - 1/k + o(1)){n \choose 2}}$, as claimed.
\end{proof}

Theorems~\ref{order} and \ref{highord} determine the possible speeds of a hereditary property of ordered graphs below $2^{n-1}$ and above $2^{cn^2}$, but in the large range in between many questions remain. In \cite{BBMpar} the current authors conjectured that for such a property $\P$ either $|\P_n| < c^n$ for some constant $c$, or $|\P_n| \ge \sum_{k=0}^{\lfloor n/2 \rfloor} {n \choose {2k}} k!$ for every $n \in \N$. They also proved several special cases of the conjecture (some of which were proved independently by Klazar and Marcus~\cite{KM}), each of which generalizes the well-known Stanley-Wilf conjecture, which was recently proved by the combined results of Klazar~\cite{Klaz1} and Marcus and Tardos~\cite{MT}.

We shall discuss these and other open questions in greater length in Section~\ref{probsec}, but let us now return to the proof of Theorem~\ref{order}. The proof will require some notation, and for convenience we shall give a portion of it here, before we begin.

Let $G$ be an ordered graph with $V(G) = [n]$. The \emph{length} of the edge $ij \in E(G)$ is $|i - j|$, and $G$ is \emph{$\ell$-complete} if it has all edges of length at least $\ell$, and \emph{$\ell$-empty} if it has none. If $x \in [n]$, and $\ell \in \N$, then let $N_\ell(x) = [x - \ell + 1, x + \ell - 1]$ (where $[a,b] = \{k \in \N : a \le k \le b\}$). Let $\Gamma(x)$ denote the set of neighbours of $x$ in $G$. We say $x$ and $y$ are \emph{$\ell$-homogeneous} (and write $x \sim_\ell y$) if $\Gamma(x) \setminus (N_\ell(x) \cup N_\ell(y)) = \Gamma(y) \setminus (N_\ell(x) \cup N_\ell(y))$, and say that $B \subset V(G)$ is an \emph{$\ell$-homogeneous block} if it is a set of consecutive vertices such that $x \sim_\ell y$ for every $x, y \in B$. Note that if $B$ is a $\ell$-homogeneous block then $G[B]$ is $\ell$-complete or $\ell$-empty. If $B$ is a maximal $1$-homogeneous block, we say simply that it is a \emph{homogeneous block}. Note that $\sim_1$ is an equivalence relation, so the homogeneous blocks of $G$ are determined uniquely.

Let $A,B \subset [n]$. We shall write $A < B$ if $a < b$ for every $a \in A$ and $b \in B$, and $A \sim_\ell B$ if $a \sim_\ell b$ for every $a \in A$ and $b \in B$. If $G$ is an ordered graph, and $A, B \subset V(G)$ with $A < B$, let $G[A]$ denote the ordered graph induced by the set $A$, and let $G[A,B]$ denote the bipartite ordered graph induced by the edges between $A$ and $B$. We shall write $G - A$ for $G[[n] \setminus A]$, and $G[a,b]$ for $G[[a,b]]$ if $a,b \in \N$. Finally, $G \le H$ will mean (as above) that $G$ is an induced ordered subgraph of $H$.

The rest of the paper is organised as follows. In Section~\ref{keysec} we shall prove the key lemma in the proof of Theorem~\ref{order}; in Section~\ref{typesec} we shall prove that the existence of certain structures in $\P$ implies that the speed is at least $2^{n-1}$; in Section~\ref{polysec} we shall prove the jump from polynomial speed to $F_n$; in Section~\ref{ellsec} we shall prove various lemmas about $\ell$-empty ordered graphs; in Section~\ref{Fnksec} we shall deduce the structure of a property with speed $p(n)F_{n,k}$; in Section~\ref{pfsec} we shall prove Theorem~\ref{order}; and in Section~\ref{probsec} we shall discuss some further problems, including the possible exponential speeds above $2^{n-1}$.

\section{The key lemma}\label{keysec}

We start by defining, for each pair $k,\ell \in \N$, ten basic structures. The structures come in four flavours.\\[+1ex]
Type 1: there are vertices $y$ and $x_1 < \ldots < x_{2k}$ in $G$ such that $y < x_1$ or $y > x_{2k}$, and for $1 \le i < 2k$, $yx_i \in E(G)$ iff $yx_{i+1} \notin E(G)$.\\[+1ex]
Type $2(a)$: there are vertices $x_1 < \ldots < x_{2k} < y_1 < \ldots < y_{2k}$ in $G$ such that $x_iy_i \in E(G)$ iff $x_{i+1}y_{i+1} \notin E(G)$.\\[+1ex]
Type $2(b)$: there are vertices $x_1 < \ldots < x_{2k} < y_{2k} < \ldots < y_1$ in $G$ such that $x_iy_i \in E(G)$ iff $x_{i+1}y_{i+1} \notin E(G)$.\\[+1ex]
Type 3: there are vertices $x_1 < z_{1,1} < \ldots < z_{1,\ell - 1} < y_1 < x_2 < z_{2,1} < \ldots < z_{2,\ell-1} < y_2 < x_3 < \ldots < y_{2k-1} < x_{2k} < z_{2k,1} < \ldots < z_{2k,\ell-1} < y_{2k}$ in $G$ such that $x_iy_i \in E(G)$ iff $x_{i+1}y_{i+1} \notin E(G)$.\\

Note that there are four different structures of Type 1, and two each of Types 2(a), 2(b) and 3. We refer to these as $k$-structures of Type 1 and 2, and $(k,\ell)$-structures of Type 3 (throughout we shall say ``Type 2" when we mean ``Type 2(a) or Type 2(b)"). They are not graphs, but sub-structures contained in graphs: instead of saying that ``a structure of Type $i$ occurs in $\P$'' it would be more precise to say that ``there is a graph $G \in \P$ admitting a structure of Type $i$". However, for smoothness of presentation, we sometimes handle them as graphs.

The key lemma in the proof of Theorem~\ref{order} will be the following.

\begin{lemma}\label{key}
Let $k,\ell \in \N$, and $G$ be any ordered graph. If $G$ contains no $k$-structure of Type 1 or 2, and no $(k,\ell)$-structure of Type 3, then the vertices of $G$ may be partitioned into blocks $B_1 < \ldots < B_m$, with $m \le 256k^4$, and each block $\ell$-homogeneous.
\end{lemma}

\begin{proof}
Let $k,\ell \in \N$, and $G$ be an ordered graph with vertex set $[n]$. Suppose that no $k$-structure of Type 1 or 2, and no $(k,\ell)$-structure of Type 3 occurs in $G$, and suppose without loss of generality that $(1,\ell+1) \in E(G)$. Let $i_0 = 0$, $i_1 = \ell + 1$, and let $i_2$ be minimal under the condition that it is an endpoint of a non-edge $ji_2$, with $i_1 = \ell + 1 < j \le i_2 - \ell$. Now fixing $i_2$, let $i_3$ be the minimal number under the condition that it is an endpoint of an edge $ji_3$, with $i_2 < j \le i_3 - \ell$. If no such edge / non-edge exists at stage $t$, then set $i_t = n + 1$ and stop. Continuing in this way, the sequence $\{i_0, \ldots, i_t\}$ may be defined, but it may not continue further than $t = 2k$, otherwise a $(k,\ell)$-structure of Type 3 would appear as a sub-structure. For each $j \in [t]$, the graph spanned by $[i_{j - 1}+1, i_j-1]$ is either $\ell$-complete or $\ell$-empty, depending on the parity of $j$. This means that the vertex set of $G$ can be partitioned into at most $2k-1$ vertices and at most $2k$ blocks $A_1 < \ldots < A_t$ of consecutive vertices, where the blocks span $\ell$-complete or $\ell$-empty graphs.

Let $j \in [t]$ and consider the block $A_j$. Let $A_j = [u_j,v_j]$, and let $S_j = \{ s \in [u_j, v_j-1] : \Gamma(s) \setminus A_j \neq \Gamma(s+1) \setminus A_j \}$ be the set of vertices `separated' from the next vertex to the right by a vertex outside $A_j$. We shall show that $|S_j| < 64k^3$.

For each $s \in S_j$, choose a vertex $w = w(s) \in [n] \setminus A_j$ such that $sw \in E(G)$ but $(s+1)w \notin E(G)$, or vice versa. Let $T_j = \{w(s) : s \in S_j\}$. If any vertex $w$ is chosen by more than $2k$ different vertices of $S_j$ then $G$ contains a $k$-structure of Type 1, a contradiction. Hence if $|S_j| \ge 64k^3$, then $|T_j| \ge 32k^2$. At least $16k^2$ of these vertices must lie to the left, say, of $A_j$. Denote these vertices $w_1 < \ldots < w_{16k^2}$, and let $f: [16k^2] \to S_j$ map $x \in [16k^2]$ to a vertex which chose $w_x$, i.e., $w(f(x)) = w_x$ for each $x \in [16k^2]$. By the Erd\H{o}s-Szekeres Theorem, there is a subsequence $A$ of this set with length at least $4k$ on which $f$ is monotone. Note that $f$ is injective, so $f$ is in fact strictly monotone on $A$. If $f$ is increasing on $A$, then $G$ admits a $k$-structure of Type $2(a)$; if it is decreasing then $G$ admits a $k$-structure of Type $2(b)$. In either case we contradict one of our assumptions, so $|S_j| < 64k^3$ as claimed.

Now, partition $[n]$ into sets $\{B_1, \ldots, B_{2m + 1}\}$ of consecutive vertices, with $m < 128k^4$, as follows. Let $\bigcup_{j = 1}^t S_j \cup \{i_1, \ldots, i_{t-1}\}$ have elements $a_1 < \ldots < a_m$. Note that $m \le (64k^3 - 1)t + t - 1 \le 128k^4 - 1$ by the comments above, and let $a_0 = 0$ and $a_{m+1} = n + 1$. For each $j \in [m]$, let $B_{2j} = a_j$, and for each $j \in [m+1]$, let $B_{2j-1} = [a_{j - 1} + 1, a_j - 1]$.

Now, each set $B_i$ either consists of a single vertex, or $B_i \subset A_j \setminus S_j$ for some $j \in [t]$, and is an interval of $[n]$. Hence if $x, y \in B_i$, then $\Gamma(x) \setminus A_j = \Gamma(y) \setminus A_j$. Since $A_j$ is $\ell$-complete or $\ell$-empty, it follows that $x \sim_\ell y$. So we have partitioned $V(G)$ into at most $2m + 1 < 256k^4$ blocks of consecutive vertices, with each block $\ell$-homogeneous, as required.
\end{proof}

\section{Structures of Type 1, 2 and 3}\label{typesec}

In order to use Lemma~\ref{key} to prove Theorem~\ref{order}, we must give sharp lower bounds on the possible speeds of hereditary properties containing large structures of Type 1, 2 or 3. The bounds are given by Lemmas~\ref{t1t2} and \ref{t3}. Lemma~\ref{t1t2} will show that if $\P$ contains arbitrarily large structures of Type 1 or 2, then $|\P_n| \ge 2^{n-1}$ for every $n \in \N$. To prove it we will need to handle some particular ordered graphs, and for ease of presentation we first define them here.

Let $m \in \N$, let $X = \{x_1, \dots, x_m\}$ and $Y = \{y_1, \dots, y_m\}$ be disjoint sets of vertices satisfying $x_1 < \ldots < x_m < y_1 < \ldots < y_m$, and let $I = (I_1,I_2,I_3,I_4) \in \{0,1\}^4$. The graph $M^<_I = M^<_I(X,Y)$ on $X \cup Y$ has the edge set defined as follows:
\begin{enumerate}
\item[$(i)$] $x_ix_j \in E(M^<_I)$ if and only if $I_1 = 1$;
\item[$(ii)$] if $i < j$, then $x_iy_j \in E(M^<_I)$ if and only if $I_2 = 1$;
\item[$(iii)$] if $i > j$, then $x_iy_j \in E(M^<_I)$ if and only if $I_3 = 1$;
\item[$(iv)$] $y_iy_j \in E(M^<_I)$ if and only if $I_4 = 1$;
\item[$(v)$] $x_iy_i \in E(M^<_I)$ if and only if $i$ is odd.
\end{enumerate}
The graph $M^>_I = M^>_I(X,Y)$ is obtained in exactly the same way, except the vertices of $Y$ are first renamed so that $y_1 > \ldots > y_m$.

\begin{lemma}\label{justcount}
Let $n \in \N$, $|X| = |Y| \ge n^2 + n$ and $I \in \{0,1\}^4$. Then $M^<_I(X,Y)$ and $M^>_I(X,Y)$ each have at least $2^{n-1}$ distinct induced ordered subgraphs of order $n$.
\end{lemma}

\begin{proof}
Let $I = (I_1,I_2,I_3,I_4)$ and $M = M^<_I(X,Y)$. Assume first, by taking the complementary graph (and removing $x_1$ and $y_1$) if necessary, that $I_2 = I_3 = 1$ does not hold, and that either $I_2 = I_3 = 0$ or $I_1 = 0$. In the latter case, we may assume also (by symmetry) that $I_3 = 0$. Note first that the result is clear if $n \le 2$. Now, if $n = 3$ and $I_2 = I_3 = I_4 = 0$, then the four ordered subgraphs induced by the sets $\{x_1,y_1,y_2\}$, $\{x_1,x_2,y_1\}$, $\{x_1,x_3,y_3\}$ and $\{y_1,y_2,y_3\}$ are distinct, and if $n = 3$, $I_2 = I_3 = 0$ and $I_4 = 1$, then the four ordered subgraphs induced by the sets $\{x_1,y_1,y_2\}$, $\{x_1,y_2,y_3\}$, $\{x_1,x_2,y_1\}$ and $\{y_1,y_2,y_3\}$ are distinct, so we are done in these cases as well. Hence we may assume that either $n \ge 4$ and $I_2 = I_3 = 0$, or $n \ge 3$, $I_2 = 1$ and $I_1 = I_3 = 0$.

We shall describe an injective map $\phi$ from the subsets of $[n]$ of even size to induced subgraphs of $M$ on $n$ vertices; since there are $2^{n-1}$ such subsets, this will suffice to prove the lemma. Given a subset $S = \{s_1, \ldots, s_\ell, t_1, \ldots, t_\ell\} \subset [n]$, with $1 \le s_1 < \ldots s_\ell < t_1 < \ldots < t_\ell \le n$, we shall define $\phi(S)$ to be a subgraph $G$ of $M$ with a matching, or `star-matching', between the vertices $\{s_1, \ldots, s_\ell\}$ and $\{t_1, \ldots, t_\ell\}$ (in a star-matching the edge-set is $\{s_it_j : i \ge j\}$ or $\{s_it_j : i \le j\}$). This will allow us to reconstruct $S$ from $G$, and hence show that $\phi$ is injective.

To be precise, first let
$$A = \{x_i, y_i : i = 2jn - 1, 1 \le j \le \ell \},$$
so $x_iy_i \in E(M)$ for each $i$ with $x_i, y_i \in A$. The vertices of $A$ will correspond to the elements of $S$. We now need to fill in the space between the vertices of $A$, but without creating any edges that will prevent us from identifying $S$ (see Figure 1). To this end, let $s_0 = 0$, $s_{\ell+1} = t_1$ and $t_{\ell+1} = n+1$, and let
\begin{align*}
B & = \; \{x_i : i \in [(2j-1)n + 1, (2j-1)n + s_j - s_{j - 1} - 1]\textup{ for some }j \in [\ell + 1]\} \\
& \; \cup \; \{y_i : i \in [2jn + 1, 2jn + t_{j + 1} - t_j - 1]\textup{ for some }j \in [\ell] \}.
\end{align*}
This is possible because $|X|, |Y| \ge n^2 + n - 1$.\\

\[ \unit = 0.45cm \hskip -18\unit
\medline \dl{0}{2}{1}{2} \dl{1}{2}{1}{0} \dl{0}{0}{1}{0} \dl{0}{0}{0}{2}
\thnline \dl{0.1}{1.9}{0.9}{0.1} \dl{0.1}{1}{0.5}{0.1} \dl{0.5}{1.9}{0.9}{1}
\varpt{5000} \pt{0.5}{-1} \point{-1.6}{-1.8}{$x_{2n-1}$} \medline \bez{0.5}{-1}{6.5}{1}{12.5}{-1} \pt{12.5}{-1} \point{13}{-0.4}{$y_{2n-1}$}
\point{12}{-4}{$ \medline \dl{0}{2}{1}{2} \dl{1}{2}{1}{0} \dl{0}{0}{1}{0} \dl{0}{0}{0}{2}
\thnline \dl{0.1}{1.9}{0.9}{0.1} \dl{0.1}{1}{0.5}{0.1} \dl{0.5}{1.9}{0.9}{1} $}
\point{0}{-7}{$ \medline \dl{0}{2}{1}{2} \dl{1}{2}{1}{0} \dl{0}{0}{1}{0} \dl{0}{0}{0}{2}
\thnline \dl{0.1}{1.9}{0.9}{0.1} \dl{0.1}{1}{0.5}{0.1} \dl{0.5}{1.9}{0.9}{1}
\varpt{5000} \pt{0.5}{-1} \point{-1.6}{-1.8}{$x_{4n-1}$} \medline \bez{0.5}{-1}{6.5}{1}{12.5}{-1} \pt{12.5}{-1} \point{13}{-0.4}{$y_{4n-1}$}
\point{12}{-4}{$ \medline \dl{0}{2}{1}{2} \dl{1}{2}{1}{0} \dl{0}{0}{1}{0} \dl{0}{0}{0}{2}
\thnline \dl{0.1}{1.9}{0.9}{0.1} \dl{0.1}{1}{0.5}{0.1} \dl{0.5}{1.9}{0.9}{1} $} $}
\thkline \dline{6.5}{-11.5}{6.5}{-14.5}{4}
\point{0}{-17}{$ \medline \dl{0}{2}{1}{2} \dl{1}{2}{1}{0} \dl{0}{0}{1}{0} \dl{0}{0}{0}{2}
\thnline \dl{0.1}{1.9}{0.9}{0.1} \dl{0.1}{1}{0.5}{0.1} \dl{0.5}{1.9}{0.9}{1}
\varpt{5000} \pt{0.5}{-1} \point{-1.6}{-1.8}{$x_{2\ell n-1}$} \medline \bez{0.5}{-1}{6.5}{1}{12.5}{-1} \pt{12.5}{-1} \point{13}{-0.4}{$y_{2\ell n-1}$}
\point{12}{-4}{$ \medline \dl{0}{2}{1}{2} \dl{1}{2}{1}{0} \dl{0}{0}{1}{0} \dl{0}{0}{0}{2}
\thnline \dl{0.1}{1.9}{0.9}{0.1} \dl{0.1}{1}{0.5}{0.1} \dl{0.5}{1.9}{0.9}{1} $} $}
\point{0}{-24}{$ \medline \dl{0}{2}{1}{2} \dl{1}{2}{1}{0} \dl{0}{0}{1}{0} \dl{0}{0}{0}{2}
\thnline \dl{0.1}{1.9}{0.9}{0.1} \dl{0.1}{1}{0.5}{0.1} \dl{0.5}{1.9}{0.9}{1} $}
\point{17}{-13}{$ \medline \dl{0}{2}{1}{2} \dl{1}{2}{1}{0} \dl{0}{0}{1}{0} \dl{0}{0}{0}{2}
\thnline \dl{0.1}{1.9}{0.9}{0.1} \dl{0.1}{1}{0.5}{0.1} \dl{0.5}{1.9}{0.9}{1} \point{1.1}{0.8}{ $= B$} $}
\point{2}{-27}{Figure 1: The set $A \cup B$}
\]\\

Define $\phi(S) = M[A \cup B]$. Notice that this gives $\phi(\emptyset) = E_n$ or $K_n$, depending on whether $I_1 = 0$ or $1$. Also, if $S \neq \emptyset$ then $s_j - s_{j - 1} \le n - 1$ for every $j \in [\ell+1]$, and $t_{j + 1} - t_j  \le n - 1$ for every $j \in [\ell]$, so $A \cap B = \emptyset$. Therefore
\begin{align*}
|A \cup B| & = \; 2\ell + \sum_{j = 1}^{\ell + 1} (s_j - s_{j - 1} - 1) + \sum_{j = 1}^\ell (t_{j + 1} - t_j - 1)\\
& = \; 2\ell + (s_{\ell + 1} - \ell - 1) + (t_{\ell + 1} - t_1 - \ell) \; = \; n.
\end{align*}
Moreover, if we identify $A \cup B$ with $[n]$ in the obvious way, then $A = S$. We have two cases to investigate.\\[+2ex]
\ul{Case 1}: $I_2 = I_3 = 0$.\\[+2ex]
Let $G = \phi(S)$ for some even-size subset $S$ of $[n]$, so $G$ is an ordered graph with vertex set $[n]$. We wish to show that $S$ is uniquely determined by $G$. Let $S = \{s_1, \ldots, s_\ell, t_1, \ldots, t_\ell\}$ as before, and recall from above that we may assume that $n \ge 4$.

If $G = K_n$ then the vertices of $G$ all came from $X$, or all from $Y$, since $I_2 = I_3 = 0$ and $n \ge 3$. Thus in this case $S = \emptyset$. Also, if $|S| \ge 2$ then $s_1t_1 \in E(G)$, so if $G = E_n$ then we also have $S = \emptyset$ (there is no contradiction here -- only one of the two cases is possible for a given $M$). Therefore we are done in the case $G \in \{K_n, E_n\}$, so assume that $G$ has at least one edge and one non-edge, and hence that $S \neq \emptyset$.

Suppose first that $\Gamma_G(1) \neq \{2\}$, and recall that $G[1,t_1 - 1]$ must be complete or empty. We claim that $t_1$ must be the left-most vertex of $G$ with a neighbour to its left, but which is not part of a clique involving all of the vertices to its left. To see this, observe that $t_1$ certainly has a neighbour to its left, since $s_1t_1 \in E(G)$, and that $s_1$ is its only neighbour to its left, since $I_2 = I_3 = 0$. Hence if $t_1$ is part of a clique involving all the vertices to its left, then $t_1 = 2$, so $S = \{1,2\}$ and $\Gamma_G(1) = \{2\}$, a contradiction. Now suppose some vertex to the left of $t_1$ has a neighbour to its left. Then $I_1 = 1$, so $G[1,t_1 - 1]$ is a clique. This proves the claim, so if $\Gamma_G(1) \neq \{2\}$ then we can reconstruct $t_1$. But now $S = \{u,v : uv \in E(G)$ and $u < t_1 \le v \}$ is the only possibility for $S$, since $I_2 = I_3 = 0$, so the only edges of $\phi(S)$ between $X$ and $Y$ are the edges $s_it_i$, for $i \in [\ell]$.

So suppose next that $\Gamma_G(1) = \{2\}$, and suppose also that $\Gamma_G(2) \neq \{1,3\}$. We shall prove that in this case $S = \{1,2\}$. Indeed, since $G[1,t_1 - 1]$ must be complete or empty, we have $t_1 \in \{2,3\}$. Now, if $t_1 = 3$ then $s_1 = 2$, since $s_1t_1 \in E(G)$ and $13 \notin E(G)$. But then $\Gamma_G(2) = \{1,3\}$, since $I_2 = I_3 = 0$, which is a contradiction. Hence $t_1 = 2$, and so $S = \{1,2\}$.

Suppose finally that $\Gamma_G(1) = \{2\}$ and $\Gamma_G(2) = \{1,3\}$. Since $G[1,t_1 - 1]$ must be complete or empty, we have $t_1 \in \{2,3\}$, and since $G[t_1,n]$ must also be complete or empty and $n \ge 4$, $t_1 = 2$ is impossible. Thus $t_1 = 3$, and so $s_1 \neq 1$, since $3 \notin \Gamma_G(1)$. Hence $S = \{2,3\}$ in this case.\\

By the comments above, we have reduced the problem to the following case.\\[+2ex]
\ul{Case 2}: $I_1 = 0$, $I_2 = 1$ and $I_3 = 0$.\\[+2ex]
Recall that we may assume $n \ge 3$. Observe that $1t_1 \in E(G)$, so $\Gamma_G(1) = \emptyset$ if and only if $S = \emptyset$. Therefore we may assume that $S \neq \emptyset$.

Consider the homogeneous blocks $B_1, \ldots, B_k$ of $G$. We claim that they are exactly the sets $\{1, \ldots, s_1\}, \{s_1 + 1, \ldots, s_2\}, \ldots, \{s_\ell + 1 , \ldots, t_1 - 1\}, \{t_1, \ldots, t_2 - 1\}, \ldots, \{t_\ell, \ldots, n\}$. To see this, consider any two vertices $u,v \in G$ and consider the following cases.
\begin{enumerate}
\item[$(i)$] $s_{j-1} < u < v \le s_j$ for some $j \in [\ell]$. Then $\Gamma(u) \setminus \{v\} = \Gamma(v) \setminus \{u\} = [t_j,n]$, so $u \sim v$.
\item[$(ii)$] $s_\ell < u < v < t_1$. Then $\Gamma(u) = \Gamma(v) = \emptyset$, since $I_1 = I_3 = 0$, and $u$ and $v$ are `below' every vertex of $B \cap Y$. So $u \sim v$.
\item[$(iii)$] $t_j \le u < v < t_{j+1}$ for some $j \in [\ell]$ and $I_4 = 0$. Then $\Gamma(u) \setminus \{v\} = \Gamma(v) \setminus \{u\} = [1,s_j]$, so $u \sim v$.
\item[$(iv)$] $t_j \le u < v < t_{j+1}$ for some $j \in [\ell]$ and $I_4 = 1$. Then $\Gamma(u) \setminus \{v\} = \Gamma(v) \setminus \{u\} = [1,s_j] \cup [t_1,n] \setminus \{u,v\}$ so $u \sim v$.
\item[$(v$)] $u \le s_j < v < t_1$ for some $j \in [\ell]$. Then $ut_j \in E(G)$ but $vt_j \notin E(G)$, so $u \not\sim v$.
\item[$(vi)$] $t_1 \le u < t_j \le v$. Then $us_j \notin E(G)$ but $vs_j \in E(G)$, so $u \not\sim v$.
\item[$(vii)$] $u = t_1 - 1$ and $v = t_1$. Then $1u \notin E(G)$ and $1v \in E(G)$, so $u \not\sim v$.
\end{enumerate}

Now, there are either $2\ell$ or $2\ell + 1$ homogeneous blocks in $G$ (depending on whether or not $t_1 = s_\ell + 1$), and in either case the set $S$ must consist of the right-most vertices of the first $\ell$ blocks and the left-most vertices of the last $\ell$. Thus we have proved that in all cases we can reconstruct $S$ from $\phi(S)$, and hence $\phi$ is injective. This proves that $M$ has at least $2^{n-1}$ distinct induced subgraphs, as claimed.

The proof for $M^>_I(X,Y)$ is almost identical.
\end{proof}

\noindent Lemma~\ref{justcount} and Ramsey's Theorem now give the following result.

\begin{lemma}\label{t1t2}
Let $\P$ be a hereditary property of ordered graphs. Suppose a $k$-structure of Type 1 or 2 occurs in $\P$ for arbitrarily large values of $k$. Then $|\P_n| \ge 2^{n-1}$ for every $n \in \N$.
\end{lemma}

\begin{proof}
Let $\P$ be a hereditary property of ordered graphs, and suppose first that $\P$ contains a $k$-structure of Type 1 for arbitrarily large $k$. Let $n \in \N$ and choose $k \ge n - 1$ and an ordered graph $G \in \P$ containing a $k$-structure of Type 1 on vertices $\{y, x_1, \ldots, x_{2k}\}$. Without loss of generality, assume that $y < x_1 < \ldots < x_{2k}$, and $yx_1 \in E(G)$. Now, for each subset $S \subset [n-1]$, let the ordered graph $G_S \in \P_n$ be induced by the vertices $y \cup \{x_i : i \in (2S - 1) \cup (2S^c)\}$, where $2S - 1 = \{2s - 1 : s \in S\}$ and $2S^c = \{2s : s \in [n - 1] \setminus S\}$. The graphs $G_S$ are all distinct, since the set $S$ can be recovered from $G_S$ by considering the neighbours of the left-most vertex, and are all in $\P_n$. It follows that $|\P_n| \ge 2^{n-1}$.

So now assume that for some $K \in \N$ there is no $k$-structure of Type 1 in $\P$ for $k \ge K$, and $\P$ contains a $k$-structure of Type 2 for arbitrarily large $k$. It can be shown fairly easily that an $n$-structure of Type 2 contains at least $2^{n/2}$ distinct ordered subgraphs on $n$ vertices. To do better than this we will use Ramsey's Theorem to produce some uniformity on the unknown edges. Let $R_r(s)$ be the smallest number $m$ such
that any $r$-colouring of the edges of $K_m$ contains a monochromatic $K_s$. Let $n \in \N$, $r = 2^{16}$, $k = R_r(\max\{n^2+n, K+1\})$, and choose an ordered graph $G \in \P$ containing a $k$-structure of Type 2 on vertices $\{x_1, \ldots, x_{2k}, y_1, \ldots, y_{2k}\}$. We assume without loss of generality that $x_1 < \ldots < x_{2k} < y_i$ for every $i \in [2k]$, that \emph{either} $y_1 < \ldots < y_{2k}$ \emph{or} $y_1 > \ldots > y_{2k}$, and that $x_iy_i \in E(G)$ if and only if $i$ is odd. We shall thus prove the result for both Type 2(a) and Type 2(b) properties at the same time. Let $X = \{x_1, \ldots, x_{2k}\}$ and $Y = \{y_1, \ldots, y_{2k}\}$.

We first split our Type 2 structure up into blocks of four vertices each, as follows: $D_1 = \{x_1, y_1, x_2, y_2\}, D_2 = \{x_3, y_3, x_4, y_4\} , \ldots, D_k = \{x_{2k-1}, y_{2k-1}, x_{2k}, y_{2k}\}$, and let $J$ be the complete graph with these $k$ blocks as vertices. Define a $2^{16}$-colouring on the edges of $J$ by associating the bipartite ordered graph $G[D_i,D_j]$ with the edge $D_iD_j$. By our choice of $k$, there must be a complete monochromatic subgraph of $J$ on $s \ge \max\{n^2+n, K+1\}$ blocks. By renaming the vertices of $G$ if necessary, we may assume that these blocks are $D_1, \ldots, D_s$.

Each pair of blocks $D_i, D_j$ (with $1 \le i < j \le s$) induces the same bipartite ordered graph $G[D_i,D_j]$. The next claim shows that there are only a small number of possibilities for this graph.\\

\noindent\ul{Claim 1}: Each of the ordered bipartite graphs $G[\{x_1,x_2\}, \{x_3,x_4\}]$, $G[\{x_1,x_2\}, \{y_3,y_4\}]$, $G[\{x_3,x_4\}, \{y_1,y_2\}]$ and $G[\{y_1,y_2\}, \{y_3,y_4\}]$ is either complete or empty.

\begin{proof}
Let $u \in D_1$, and suppose that $\Gamma(u) \cap D_2 \notin \{\emptyset, D_2 \cap X, D_2 \cap Y, D_2\}$. We shall show that $\P$ contains a $K$-structure of Type 1.

Let $\{v,w\} = D_2 \cap X$ with $v < w$, and suppose that $uv \in E(G)$ and $uw \notin E(G)$ (the proof in the other cases is identical). Consider the set $\{u, x_3, x_4, \ldots, x_{2s}\}$, and note that either $u < x_3$ or $u > x_{2s}$. Now, since $G[D_i,D_j]$ is the same for each $1 \le i < j \le s$, we have $ux_i \in E(G)$ if and only if $i$ is odd. Also, $2s \ge 2K + 2$. Thus the graph $G[\{u, x_3, x_4, \ldots, x_{2s}\}] \in \P$ contains a $K$-structure of Type 1.

We now have a contradiction, so in fact $\Gamma(u) \cap D_2 \in \{\emptyset, D_2 \cap X, D_2 \cap Y, D_2\}$ for every $u \in D_1$. Similarly, one can show that $\Gamma(u) \cap D_1 \in \{\emptyset, D_1 \cap X, D_1 \cap Y, D_1\}$ for every $u \in D_2$. The result now follows easily.
\end{proof}

Now, let $X' = \{x_i : i \in [2s]$ and $i \equiv 1,4 \pmod 4\}$ and $Y' = \{y_i : i \in [2s]$ and $i \equiv 1,4 \pmod 4\}$, and let $H$ be the subgraph of $G$ induced by the vertices of $X' \cup Y'$. Note that no two distinct vertices of $X'$ lie in the same block $D_i$, and similarly for $Y'$.\\

\noindent\ul{Claim 2}: For some $I \in \{0,1\}^4$, $H = M^<_I(X',Y')$ or $M^>_I(X',Y')$.

\begin{proof}
It is clear from Claim 1, and the fact that $G[D_i,D_j]$ is the same for each $1 \le i < j \le s$, that $H[X']$ and $H[Y']$ are either complete or empty. So let $x_i \in X'$ and $y_j \in Y'$, and observe that
\begin{itemize}
\item[$(i)$] if $i < j$ then $x_iy_j \in E(H)$ if and only if $G[\{x_1,x_2\}, \{y_3,y_4\}]$ is complete,
\item[$(ii)$] if $i > j$, then $x_iy_j \in E(H)$ if and only if $G[\{x_3,x_4\}, \{y_1,y_2\}]$ is complete,
\item[$(iii)$] if $i = j$, then $x_iy_j \in E(H)$ if and only if $i = j$ is odd.
\end{itemize}
Therefore $H = M^<_I(X',Y')$ or $M^>_I(X',Y')$ for some $I \in \{0,1\}^4$.
\end{proof}

We have shown that for some $I \in \{0,1\}^4$, either $M^<_I(X',Y') \in \P$ or $M^>_I(X',Y') \in \P$, with $|X'| = |Y'| \ge n^2 + n$. By Lemma~\ref{justcount}, it follows that $|\P_n| \geq 2^{n-1}$. Since $n$ was arbitrary, the result follows.
\end{proof}

For Type 3 structures, a different bound holds.

\begin{lemma}\label{t3}
Let $\ell \in \N$, and $\P$ be a hereditary property of ordered graphs. Suppose a $(k,\ell)$-structure of Type 3 occurs in $\P$ for arbitrarily large values of $k$. Then $|\P_n| \ge F_{n,\ell+1}$ for every $n \in \N$.
\end{lemma}

\begin{proof}
Let $\ell \in \N$, and let $\P$ be a hereditary property of ordered graphs containing $(k,\ell)$-structures of Type 3 for arbitrarily large values of $k$. If $\P$ also contains $k$-structures of Type 1 for arbitrarily large values of $k$, then by Lemma~\ref{t1t2}, $|\P_n| \ge 2^{n-1} \ge F_{n,\ell+1}$ for every $n \in \N$, in which case we are done. So assume that there exists $K \in \N$ such that there is no $k$-structure of Type 1 in $\P$ for any $k \ge K$.

Let $n \in \N$, $r = 2^{4(\ell+1)^2}$, $k = R_r(\max\{2n,2K+1\})$, and choose a graph $G \in \P$ containing a $(k,\ell)$-structure of Type 3. Let the vertices of this Type 3 structure be $$\{x_i : i \in [2k]\} \cup \{ y_i : i \in [2k] \} \cup \{ z_{i,j} : i \in [2k], j \in [\ell-1] \},$$ where $x_1 < z_{1,1} < \ldots < z_{1,\ell-1} < y_1 < \ldots < x_{2k} < z_{2k,1} < \ldots < z_{2k,\ell-1} < y_{2k}$, and without loss of generality $x_iy_i \in E(G)$ if and only if $i$ is odd. We shall apply the same method as in the proof of Lemma~\ref{t1t2}.

As before, group the vertices into blocks, this time of size $2(\ell+1)$, as follows: let $D_i = \{x_{2i-1}, y_{2i - 1}, x_{2i}, y_{2i}\} \cup \{z_{i,j} : i \in \{2i - 1, 2i\}, j \in [\ell - 1]\}$ for $1 \le i \le k$. Note that $D_1 < D_2 < \ldots < D_k$. Let $J$ be the complete graph with these $k$ blocks as vertices, and define a $2^{4(\ell+1)^2}$-colouring on the edges of $J$ by associating the bipartite ordered graph $G[D_i,D_j]$ with the edge $D_iD_j$. By our choice of $k$, there must be a complete monochromatic subgraph of $J$ on $s \ge \max\{2n, K+1\}$ blocks. By renaming the vertices of $G$ if necessary, we may assume that these blocks are $D_1, \ldots ,D_s$.

 Suppose that some vertex $u \in D_1$ sends an edge and a non-edge to $D_2$; say $uv_2 \in E(G)$ and $uw_2 \notin E(G)$, with $v_2,w_2 \in D_2$. Since $G[D_i,D_j]$ is the same for every $1 \le i < j \le s$, this means that $uv_j \in E(G)$, and $uw_j \notin E(G)$ for each $j \in [2,s]$, where $v_j$ and $w_j$ are the vertices of $D_j$ corresponding to $v_2$ and $w_2$ respectively. Thus $u$ sends an edge and a non-edge (in the same order) to each $D_j$ with $j \in [2,s]$, so $G[\{u, v_2, w_2, \ldots, v_s, w_s\}] \in \P$ contains an $(s-1)$-structure of Type 1. Since $s - 1 \ge K$, this is a contradiction, so either $D_2 \subset \Gamma(u)$, or $D_2 \subset \overline{\Gamma(u)}$, for each $u \in D_1$. Similarly, one can show that either $D_1 \subset \Gamma(u)$ or $D_1 \subset \overline{\Gamma(u)}$ for each $u \in D_2$.

It follows easily that the ordered bipartite graph $G[D_1,D_2]$ is complete or empty. Since the ordered graph $G[D_i,D_j]$ is the same for every $i < j$, this implies that either all or none of the edges $\{uv : u \in D_i, v \in D_j, 1 \le i < j \le s\}$ are in $E(G)$. Suppose, by taking the complement of $G$ if necessary, that none of these edges are in $E(G)$, and let $H$ be the subgraph of $G$ induced by the vertices $\{x_i, y_i, z_{i,t} : i \in [s]$ is odd, and $t \in [\ell-1] \}$. Note that $x_iy_i \in E(G)$ for each $x_i,y_i \in V(H)$.

We claim that $H$ has at least $F_{m,\ell+1}$ distinct induced ordered subgraphs on $m$ vertices for every $m \le n$. This is clear if $n = 1$, so let $n \ge 2$ and suppose the result is true for $n - 1$. Then $H[\{ x_i, y_i, z_{i,t} : i \in [3,s]$ is odd, and $t \in [\ell-1] \}]$ has at least $F_{m,\ell+1}$ distinct subgraphs on $m$ vertices for every $m \le n - 1$. It follows that, for $1 \le t \le \ell+1$, $H$ has at least $F_{n - t,\ell + 1}$ distinct subgraphs $M$ of order $n$ in which $\max\{v : x_1v \in E(M)\} = t$ (where $t = 1$ if $x_1$ is isolated). These are all distinct, so $H$ also has at least $F_{n-1,\ell+1} + F_{n-2,\ell+1} + \ldots + F_{n-(\ell+1),\ell+1} = F_{n,\ell+1}$ distinct subgraphs of order $n$, as claimed. Since $H \in \P$, and $n$ was arbitrary, the proof of the lemma is complete.
\end{proof}

The following corollary of Lemmas~\ref{key}, \ref{t1t2} and  \ref{t3} summarises what we have proved so far.

\begin{cor}\label{summ}
Let $\P$ be a hereditary property of ordered graphs. Suppose that $|\P_n| < F_{n,\ell+1}$ for some $\ell$ and $n \in \N$. Then there exists a $K \in \N$ such that every ordered graph $G \in \P$ may be partitioned into at most $K$ blocks of consecutive vertices, with each block $\ell$-homogeneous.
\end{cor}

\begin{proof}
Let $\P$ be a hereditary property of ordered graphs, let $\ell,n  \in \N$, and suppose that $|\P_n| < F_{n,\ell+1}$. First note that  $F_{n,\ell} \le 2^{n-1}$ for every $\ell,n \in \N$, so also $|\P_n|  < 2^{n-1}$. Hence, by Lemmas~\ref{t1t2} and \ref{t3}, there exists  $k \in \N$ such that $\P$ contains no $k$-structure of Type 1 or 2, and no $(k,\ell)$-structure of Type 3.

It now follows immediately from Lemma~\ref{key} that every ordered graph $G \in \P$ may be partitioned into at most $K = 256k^4$ blocks, with each block $\ell$-homogeneous.
\end{proof}

\section{Polynomial speed}\label{polysec}

Before considering the general case, we shall show that if $|\P_n| < F_n$ for some $n \in \N$, then $|\P_n|$ grows only polynomially, and moreover, for sufficiently large $n$ it is exactly a polynomial.

Recall that a set $B \subset V(G)$ is said to be a \emph{homogeneous block} if it is a maximal 1-homogeneous block, i.e., a maximal set of consecutive vertices such that for all $x,y \in B$, $\Gamma(x) \setminus \{y\} = \Gamma(y) \setminus \{x\}$. The {\it homogeneous  block sequence} of $G$ is the sequence $t_1 \ge t_2 \ge \ldots$, where $t_1, t_2, \dots$ are the orders of the homogeneous blocks of $G$. Note that the sequence is uniquely determined, but that $t_1, t_2, \ldots$ is not necessarily the order of the appearance of the blocks.

We need one more piece of notation before we begin. Let $G$ be an ordered graph, and $B_1, \ldots, B_m$ be a collection of 1-homogeneous blocks of $G$, with $B_1 < \ldots < B_m$ and $V(G) = B_1 \cup \ldots \cup B_m$. (For example, $B_1, \ldots, B_m$ could be the homogeneous blocks of $G$.) Define $G(B_1, \ldots, B_m)$ to be the ordered graph with possible loops, $H$, with vertex set $[m]$, and in which $ij \in E(H)$ if and only if $b_ic_j \in E(G)$ for some (and so every) $b_i \in B_i$ and $b_i \neq c_j \in B_j$. Note that a vertex of $H$ has a loop if and only if the corresponding block induces a non-trivial clique (i.e., a clique with at least two vertices).

Let $\P$ be a hereditary property of ordered graphs, and suppose that $|\P_n| < F_n$ for some $n \in \N$. By Corollary~\ref{summ}, there exists $k \in \N \cup \{0\}$ such that every ordered graph $G \in \P$ has at most $k + 1$ homogeneous blocks. Thus $t_{k+2} = 0$ for every $G \in \P$. The following lemma shows that in this case, the speed is $O(n^k)$.

\begin{lemma}\label{blocks}
Let $\P$ be a hereditary property of ordered graphs, and let $k,M \ge 0$ be integers. Suppose that for every $G \in \P$, the homogeneous block sequence of $G$ satisfies $\ds\sum_{i = k+2}^\infty t_i \le M$. Then  $|\P_n| = O(n^k)$.
\end{lemma}

\begin{proof}
Let $\P$ be a hereditary property of ordered graphs, let $k,M \ge 0$ be integers, and suppose that $t_{k+2} + t_{k+3} + \ldots \le M$ for every $G \in \P$. We shall give an upper bound on the number of ordered graphs of order $n$ in the property.

Indeed, every ordered graph $G \in \P_n$ is determined by a sequence $S = (a_1, \ldots, a_m)$ of positive integers, with $1 \le m \le k + M+1$, $\sum_{i=1}^m a_i = n$, and $\sum_{i \in I} a_i \ge n - M$ for some set $I \subset [m]$ with $|I| \le k + 1$; and an ordered graph $H$, with possible loops, on $m$ vertices. To see this, let $G \in \P_n$ have homogeneous blocks $B_1, \ldots, B_m$ satisfying $B_i < B_j$ if $i < j$, let $a_i = |B_i|$ for each $i \in [m]$, and let $H = G(B_1, \ldots, B_m)$. Now, $1 \le m \le k+M+1$, since $\sum_{i = k+2}^\infty t_i \le M$; $ \sum_{i=1}^m a_i = n$ since $|G| = n$; and $\sum_{i \in I} a_i \ge n - M$ if $I = \{i : B_i$ is one of the largest $k+1$ homogeneous blocks of $G\}$. Thus $S = (a_1, \ldots, a_m)$ and $H$ satisfy the conditions above. It is clear that $G$ can be reconstructed from $S$ and $H$.

It remains to count the number of such pairs $(S,H)$. If $m \le k+1$ then the number of sequences is just ${{n - 1} \choose {m-1}}$. If $m > k + 1$, then a sequence is determined by choosing a subset $I \subset [m]$ of size $k+1$, choosing values $\{a_i' : i \in I\}$ so that $\sum_{i \in I} a_i' = n - M$, and then partitioning $[M]$ into $m$ (possibly empty) intervals $C_1, \ldots, C_m$, and setting $a_i = |C_i|$ if $i \notin I$, and $a_i = a_i' + |C_i|$ if $i \in I$. Thus the number of sequences $S$ is at most
\begin{align*}
& \sum_{m=1}^{k+1} {{n-1} \choose {m-1}} + \sum_{m=k+2}^{k+M+1} \left( {m \choose {k+1}} {{n-M+k} \choose k} {{M + m - 1} \choose {m-1}} \right) \\
& \hspace{2.5cm} < \; (k+M+1) {{k+M+1} \choose {k+1}}  2^{2M+k} {{n+k} \choose k} \; = \; O(n^k).
\end{align*}
The number of ordered graphs $H$ with possible loops on $m$ vertices is just a constant, so this proves the result.
\end{proof}

We next show that in fact, if $k$ is taken to be minimal in Lemma~\ref{blocks}, then $|\P_n| = \Theta(n^k)$. The following lemma gives the lower bound required to prove this result. If $G$ is an ordered graph, and $u,v,w \in V(G)$, then say that $u$ and $v$ \emph{differ with respect to} $w$ if $uw \in E(G)$ but $vw \notin E(G)$, or vice-versa. This definition can be extended to homogeneous blocks in the obvious way.

\begin{lemma}\label{geton}
Let $\P$ be a hereditary property of ordered graphs, and suppose that there are ordered graphs $G \in \P$ such that $t_{k+1}$, the size of the $(k+1)^{st}$ largest homogeneous block in $G$, is arbitrarily large. Then
$$|\P_n| \ge {{n-3k-2} \choose k} = n^k/k! + O(n^{k-1})$$ as $n \to \infty$, and in particular if $k = 1$, then $|\P_n| \ge n$ for every $n \in \N$.
\end{lemma}

\begin{proof}
Let $\P$ be a hereditary property of ordered graphs, let $n,k \in \N$, and let $G \in \P$ have $k+1$ homogeneous blocks of order at least $n$. We shall construct a subgraph $H$ of $G$ with at least ${{n-3k-2} \choose k}$ distinct ordered subgraphs of order $n$. The idea is simply that $H$ should also have $k+1$ large homogeneous blocks, but at most $2k$ other vertices.

Let $B_1, \ldots, B_{k+1}$ be homogeneous blocks of $G$, each of order at least $n$, and with $B_i < B_j$ if $i < j$. Let $V_0 = B_1 \cup \ldots \cup B_{k+1}$ and let $H_0 = G[V_0]$. We shall inductively define a sequence of sets $V_0 \subset V_1 \subset \ldots \subset V_t$, for some $t \in [0,k]$, such that $|V_{i+1}| \le |V_i| + 2$ for each $i \in [1,t-1]$, and so that the sets $\{B_i : i \in [k+1]\}$ are all in different homogeneous blocks of $H = G[V_t]$.

Let $i \in [0,k-1]$, suppose we have already defined the sets $V_0 \subset \ldots \subset V_i$, and let $H_i = G[V_i]$. If the sets $\{B_i : i \in [k+1]\}$ are all in different homogeneous blocks of $H_i$, then we are done with $t = i$ and $H = H_i$. So suppose that there exists $j \in [k]$ such that $B_j$ and $B_{j+1}$ are in the same homogeneous block of $H_i$. We shall find a set $V_{i+1}$ as required, such that $B_j$ and $B_{j+1}$ are in different homogeneous blocks of $G[V_{i+1}]$. Note that $G[B_j \cup B_{j+1}]$ must be either complete or empty; without loss of generality, assume that it is empty.

Suppose first that there exists a vertex $u \in V(G) \setminus (B_j \cup B_{j+1})$ such that $B_j \subset \Gamma(u)$ but $B_{j+1} \not\subset \Gamma(u)$, or vice-versa. In this case let $V_{i+1} = V_i \cup \{u\}$. Since $B_j$ and $B_{j+1}$ differ with respect to $u$, they are in different homogeneous blocks of $G[V_{i+1}]$, as required.

So suppose that every vertex $v \in B_j \cup B_{j+1}$ has exactly the same neighbourhood in $G$. Since $B_j$ and $B_{j+1}$ are distinct homogeneous blocks of $G$, this means that there must exist vertices $v,w \in V(G)$ with $B_j < v < B_{j+1}$ and such that $v$ differs from the vertices of $B_j \cup B_{j+1}$ with respect to $w$. In this case let $V_{i+1} = V_i \cup \{v,w\}$. Again $B_j$ and $B_{j+1}$ are in distinct homogeneous blocks of $G[V_{i+1}]$, as required.

Now, the sequence $(V_0, \ldots, V_t) $ cannot continue any further than $t = k$, since if $B_j$ and $B_{j+1}$ are in different homogeneous blocks of $H_i$ (for some $i \in [0,t-1]$ and $j \in [k]$), then they are in different homogeneous blocks of $H_{i+1}$. Since each step of the process described above separates $B_j$ and $B_{j+1}$ for at least one $j \in [k]$, after $k$ steps all $k+1$ sets $B_i$ must be in different homogeneous blocks of $H = G[V_t]$.

Now, $H$ has $k+1$ homogeneous blocks of size at least $n$, and at most $2k$ other vertices, since $|V_{i+1}| \le |V_i| + 2$ for each $i \in [0,t-1]$. Consider an ordered subgraph $F$ of $H$, which includes all the vertices of $V_t \setminus V_0$, and at least two vertices from each block $B_i$. The homogeneous blocks of $F$ are $\{V(F) \cap C_i : C_i$ is a homogeneous block of $H \}$, and so two such ordered subgraphs with different sequences $(|V(F) \cap B_1|, \ldots, |V(F) \cap B_{k+1}|)$ are distinct. Hence $H$ has at least ${{n - 3k - 2} \choose k}$ distinct ordered subgraphs (this is the number of sequences of integers $(a_1, \ldots, a_{k+1})$, with $a_i \ge 2$ for each $i \in [k+1]$, and $\sum a_i = n - 2k$), and so $$|\P_n| \ge {{n - 3k - 2} \choose k} = n^k/k! + O(n^{k-1})$$ as required.

To prove the second part of the lemma, let $k = 1$ and perform the same process as above to obtain the ordered graph $H \in \P$. We are left to count the number of subgraphs of $H$ in the various different cases. Let the two large homogeneous blocks be $B$ and $C$, with $B < C$, and let $n \in \N$. There are four cases to consider.\\

\noindent Case 1: $H = G[B \cup C]$. \\

\noindent $H$ contains either all or none of the edges between $B$ and $C$; suppose without loss of generality that it contains none. Now, since $B$ and $C$ are distinct homogeneous blocks in $H$, at least one of $B$ and $C$ must induce a clique. Again without loss, suppose that $H[B]$ is complete.

For each $i \in [n]$, let $H(i)$ denote the ordered subgraph of $H$ which contains $i$ vertices from $B$, and $n - i$ vertices from $C$. The leftmost $i$ vertices of $H(i)$ induce a clique, and the leftmost $i+1$ vertices do not. Hence the ordered graphs $\{H(i) : i \in [n]\}$ are all distinct, and are all in $\P_n$. So $|\P_n| \ge n$.\\

\noindent Case 2: $H = G[B \cup C \cup \{u\}]$, where $u \notin B \cup C$, $B \subset \Gamma(u)$ and $C \not\subset \Gamma(u)$.\\

\noindent Recall that in this case (and subsequent ones) $G[B \cup C]$ is complete or empty (since $u$ was necessary to distinguish them); assume without loss that it is empty. Now consider the $n$ ordered subgraphs of $H$ obtained by taking $i$ vertices of $B$, $n - i - 1$ vertices of $C$, and $u$, for $0 \le i \le n - 1$. These subgraphs have exactly $i$ edges, so are different. So $|\P_n| \ge n$ in this case too.\\

\noindent Case 3: $H = G[B \cup C \cup \{v\}]$, where $B < v < C$, and $B \cup C \subset \Gamma(v)$.\\

\noindent Assuming again that $G[B \cup C]$ is empty, consider the $n$ ordered subgraphs of $H$ obtained by taking $i$ vertices of $B$, $n - i - 1$ vertices of $C$, and $v$, for $0 \le i \le n - 1$.  They are all distinct if $n \ge 3$, since the $(i+1)^{st}$ vertex has degree $n-1$, and all other vertices have degree 1. The result is clear if
$n \le 2$,  so in this case again $|\P_n| \ge n$ (in fact, adding the empty ordered graph, we get $|\P_n| \ge n + 1$).\\

\noindent Case 4: $H = G[B \cup C \cup \{v,w\}]$, with $B < v < C$, and $E[B \cup C \cup \{v\}] = \emptyset$.\\

\noindent If $vw \in E(H)$, then it is the only edge of $H$ (since $v$ differs from $B$ and $C$ with respect to $w$), and so the $n - 1$ ordered subgraphs of $H$ obtained by taking $i$ vertices of $B$, $n - i - 2$ vertices of $C$, and the vertices $v$ and $w$, for $0 \le i \le n-2$, are all distinct. Also none of these ordered graphs is empty, since they all contain the edge $vw$, so upon adding the empty ordered graph, we get $|\P_n| \ge n$.

Similarly, if $vw \notin E(H)$ then $B \cup C \subset \Gamma(w)$, and the same method again gives $|\P_n| \ge n$.\\

Hence $|\P_n| \ge n$ in each case. Since $n$ was arbitrary, we are done.
\end{proof}

\begin{rmk}
The constant $1/k!$ and the lower bound $n$ in Lemma~\ref{geton} are best possible. To see this, consider the family $\P$ of all ordered graphs with at most $k$ edges, each of length 1, and all independent. It is easy to see that $\P$ is hereditary, and to check that $|\P_n| = \sum_{i = 0}^k {{n-i} \choose i}$ for every $n \in \N$. We suspect that this is in fact the correct lower bound on speeds of order $n^k$.
\end{rmk}

Combining Lemmas~\ref{blocks} and \ref{geton}, we obtain the following result.

\begin{cor}\label{ntothek}
Let $\P$ be a hereditary property of ordered graphs. If $|\P_m| < F_m$ for some $m \in \N$, then $|\P_n| = \Theta(n^k)$ for some $k \in \N$, and moreover $k$ is the minimal number such that $\sum_{i = k+2}^\infty t_i$ is bounded.
\end{cor}

\begin{proof}
By Corollary~\ref{summ} with $\ell = 1$, there exists $K \in \N$ such that every ordered graph $G \in \P$ has at most $K$ homogeneous blocks. Thus $t_j = 0$ for every $j \ge K + 1$, so there exists a minimal number $k$ such that $\sum_{i=k+2}^\infty t_i$ is bounded. By Lemma~\ref{blocks}, this implies that $|\P_n| = O(n^k)$. Now since $k$ is minimal, there exist ordered graphs $G \in \P$ such that $t_{k+1}$ is arbitrarily large. Thus, by Lemma~\ref{geton}, $|\P_n| = \Omega(n^k)$, so in fact $|\P_n| = \Theta(n^k)$.
\end{proof}

We have proved that $|\P_n| = \Theta(n^k)$ for some $k \in \N$. In fact we can prove a much stronger statement, for which we will need a little preparation. We shall define a set of canonical properties, as in \cite{BBW1}, and show that if $|\P_n| = \Theta(n^k)$ with $k \in \N$, then $\P$ is the union of some subset of these properties.

Let $m \in \N$, and suppose that $H$ is an ordered graph with possible loops on $[m]$, and $b: [m] \to \N \cup \{\infty\}$ is any function. Let $\P^*(H,b)$ denote the collection of ordered graphs $G$ which may be partitioned into 1-homogeneous blocks $B_1 < \ldots < B_m$ satisfying $1 \le |B_i| \le b(i)$ for each $i$, and $G(B_1, \ldots, B_m) = H$. Define $\P(H,b)$ to be the smallest hereditary property of ordered graphs containing $\P^*(H,b)$.

Now, for each ordered graph $G$ with $t_k \neq t_{k+1}$, define the \emph{$k$-type graph} $H^{(k)}_G$ of $G$ as follows. Let $m = k + n - \sum_{i=1}^k t_i$, and partition $[n]$, the vertex set of $G$, into intervals $B_1 < \ldots < B_m$, so that either $B_i$ is one of the $k$ largest homogeneous blocks of $G$, so $|B_i| = t_j$ for some $j \in [k]$, or $|B_i| = 1$. Since $t_k \neq t_{k+1}$, the blocks $B_i$ are uniquely determined by $G$ and $k$; we shall call them the \emph{$k$-blocks} of $G$. Let $H^{(k)}_G = G(B_1, \ldots, B_m)$. Thus $H^{(k)}_G$ is uniquely determined by $G$ and $k$.

Given a set $S$ and a function $b: S \to \N \cup \{\infty\}$, let $I(b) = \{i \in S : b(i) = \infty\}$, and $J(b) = \{i \in S : b(i) > 1\}$. Let $\P$ be a property of ordered graphs, and let $G \in \P$. If $G$ has $t_k \neq t_{k+1}$, then we define the \emph{$k$-type functions} $\B_G^{(k)}$ of $G$ (with respect to $\P$) as follows. Let $\B_G^{(k)}$ be the set of functions $b : V(H^{(k)}_G) \to \N \cup \{\infty\}$ such that $J(b) = \{i : |B_i| > 1\}$, where $B_1 < \ldots < B_m$ are the $k$-blocks of $G$, $\P(H^{(k)}_G, b) \subset \P$, and $b$ is maximal subject to these constraints.

Finally, if $\P(H^{(k)}_G, b) \subset \P$ for every $b : V(H^{(k)}_G) \to \N \cup \{\infty\}$ with $J(b) = \{i : |B_i| > 1\}$, then let $b_G^{(k)}$ denote the unique function $b \in \B^{(k)}_G$. Note that $b_G^{(k)}(i) = \infty$ if $i \in J(b_G^{(k)})$, and $b_G^{(k)}(i) = 1$ otherwise.

Note that if $H^{(k)}_G = H^{(k)}_{G'}$, $b \in \B_G^{(k)}$ and $b' \in \B_{G'}^{(k)}$, then either $b = b'$, or $b$ and $b'$ are incomparable functions. We shall need the following easy observation, which is a well-known result in the theory of well-quasi orderings.

\begin{lemma}\label{func}
Let $N \in \N$, and $(b_i)_{i \in \N}$ be a sequence of functions $b_i : [N] \to \N \cup \{\infty\}$. Then the sequence contains a pair of comparable functions. In other words, there exists a pair $i,j \in \N$, with $i \neq j$, such that $b_i(n) \ge b_j(n)$ for every $n \in [N]$.
\end{lemma}

\begin{proof}
First note that, by the pigeonhole principle, we may assume that the functions are everywhere finite. We use induction on $N$. The statement is trivial for $N = 1$, so let $N \ge 2$ and assume all the functions are incomparable. Consider any function, $b_1$ say. It is not smaller than any other, so there is an index $i$ such that $b_j(i) < b_1(i)$ for infinitely many functions $b_j$, and hence there is a constant $c < b_1(i)$, such that infinitely many functions have $i^{th}$ coordinate $c$. These functions are incomparable on $[N] \setminus \{i\}$, and we are done by induction.
\end{proof}

We are ready to prove the main result of this section, that if $|\P_m| < F_m$ for some $m \in \N$, then $|\P_n|$ is exactly a polynomial for sufficiently large values of $n$.

\begin{thm}\label{poly}
Let $\P$ be a hereditary property of ordered graphs, and suppose that $|\P_m| < F_m$ for some $m \in \N$. Then there exist integers $K, N \in \N \cup \{0\}$ and $a_0, \ldots, a_K \in \Z$, such that
$$|\P_n| = \sum_{i = 0}^K a_i {n \choose i}$$ for every $n \ge N$.
\end{thm}

\begin{proof}
Let $\P$ be a hereditary property of ordered graphs, let $m \in \N$, and suppose $|\P_m| < F_m$. By Corollary~\ref{ntothek}, $|\P_n| = \Theta(n^{k - 1})$ as $n \to \infty$ for some $k \in \N$. Moreover, there exists $M \in \N$ such that $\sum_{i = k+1}^\infty t_i \le M$ for every $G \in \P$, and there exist ordered graphs $G \in \P$ with arbitrarily large values of $t_k$.

The proof is by induction on $k$. Let $k \in \N$, and assume the result holds for all smaller values of $k$. We begin by removing those ordered graphs in $\P$ for which $t_k = t_{k+1}$ is possible. Let $$\P^{(2)} = \{ G \in \P : \sum_{i = k}^\infty t_i \le 2M \},$$ and observe that $\P^{(2)}$ is hereditary. Thus $|\P^{(2)}_n| = O(n^{k-2})$, by Corollary~\ref{ntothek} (if $k = 1$ then $|\P^{(2)}_n| = 0$ for $n \ge 2M + 1$). Observe also that $t_k \ge M+1 > t_{k+1}$ for every $G \in \P \setminus \P^{(2)}$.

Next, we shall remove those ordered graphs $G \in \P \setminus \P^{(2)}$ for which $b^{(k)}_G$ is not defined, i.e., for which $|\B_G^{(k)}| \ge 2$. (Here, and throughout the proof, $\B_G^{(k)}$ and $b_G^{(k)}$ are taken with respect to $\P$.) Note that since $t_k \neq t_{k+1}$ for every $G \in \P \setminus \P^{(2)}$, $\B_G^{(k)}$ is defined for these $G$. Let $\A = \{G \in \P \setminus \P^{(2)} : |\B^{(k)}_G)| \ge 2\}$, and let
$$\P^{(3)} = \bigcup_{G \in \A, \: b \in \B_G^{(k)}} \P(H_G^{(k)}, b).$$
Now, $\P^{(3)}$ is hereditary, since $\P(H,b)$ is hereditary for every $H$ and $b$. We claim that $|\P^{(3)}_n| = O(n^{k-2})$, i.e., that $t_k$ is bounded in $\P^{(3)}$.

In order to prove the claim, let $\HH = \{H_G^{(k)} : G \in \P \setminus \P^{(2)} \}$ and observe that each $H \in \HH$ has at most $k + M$ vertices, since $|V(H_G^{(k)})| = k + n - \sum_{i = 1}^k t_i \le k + M$ for every $G \in \P \setminus \P^{(2)}$. Thus $|\HH|$ is finite. Now recall that for any $H \in \HH$, and any pair $G,G' \in \P \setminus \P^{(2)}$, if $H_G^{(k)} = H_{G'}^{(k)} = H$, $b \in \B_G^{(k)}$ and $b' \in \B_{G'}^{(k)}$, then either $b = b'$, or $b$ and $b'$ are incomparable. Thus, by Lemma~\ref{func}, there are only finitely many such functions $b$ for each $H \in \HH$, and so there are only finitely many pairs $(H,b)$ such that $H = H_G^{(k)}$ and $b \in \B_G^{(k)}$ for some $G \in \A$.

Let $\C = \{(H,b) : H = H_G^{(k)}$ and $b \in \B_G^{(k)}$ for some $G \in \A\}$, and observe that if $(H,b) \in \C$, then $|I(b)| < k$. Therefore, there is an $N(b) \in \N$ such that $t_k \le N(b)$ for every $G \in \P(H,b)$. Since $\C$ is finite, it follows that there exists an $N \in \N$ such that $t_k \le N$ for every $G \in \bigcup_{(H,b) \in \C} \P(H,b) = \P^{(3)}$, as claimed. We choose such an $N$, with $N \ge 2M + 1$. We have $\sum_{i=k}^\infty t_i \le N + M$ for every $G \in \P^{(3)}$, so let $$\P^{(4)} = \{ G \in \P : \sum_{i = k}^\infty t_i \le N + M \},$$ and observe that $\P^{(4)}$ is hereditary, and that $\P^{(2)} \cup \P^{(3)} \subset \P^{(4)}$. By Lemma~\ref{blocks}, $|\P^{(4)}_n| = O(n^{k-2})$.

We shall apply the induction hypothesis to the property $\P^{(4)}$, but first let us count the members of $(\P \setminus \P^{(4)})_n$. Let $\P^{(1)} = \P \setminus \P^{(4)}$, and consider an ordered graph $G \in \P^{(1)}$. Since $\P^{(2)} \cup \P^{(3)} \subset \P^{(4)}$, we know that $t_k > t_{k+1}$, and that $b_G^{(k)}$ is defined. Recall that $|I(b_G^{(k)})| = k$. Let $\D = \{(H,b) : H = H_G^{(k)}$ and $b = b_G^{(k)}$ for some $G \in \P^{(1)}\}$.

We claim that $G \in \P^*(H,b)$ for a unique pair $(H,b) \in \D$. Clearly $G \in \P^*(H_G^{(k)}, b_G^{(k)})$, so suppose that also $G \in \P^*(H',b')$, with $(H',b') \in \D$. Then $H' = H_{G'}^{(k)}$ and $b' = b_{G'}^{(k)}$ for some $G' \in \P^{(1)}$, and also $G(B_1, \ldots, B_m) = G'(B'_1, \ldots, B'_m) = H'$, where $B_1 < \ldots < B_m$ are 1-homogeneous blocks of $G$, with $V(G) = B_1 \cup \ldots \cup B_m$ and $1 \le |B_i| \le b'(i)$ for each $i \in [m]$, and $B'_1 < \ldots < B'_m$ are the $k$-type blocks of $G'$. Note that $|I(b')| = |J(b')| = k$.

Now, if $i \in I(b')$, then $B'_i$ is a homogeneous block of $G'$. Furthermore, if $i \in I(b')$ then $|B'_i| \ge N + 1 > 1$, since $G' \in \P^{(1)}$, and if $i \notin I(b') = J(b')$, then $|B_i| = |B'_i| = 1$. So, if $|B_i| > 1$ for each $i \in I(b')$, then (since $B'_i$ is a homogeneous block) it follows that $B_i$ must be a homogeneous block of $G$ for each $i \in I(b')$. Recall that $|H'| = m \le k + M$. We shall show that $|B_i| \ge M + 1$ for every $i \in I(b')$, so that in fact, $\{B_i : i \in I(b')\}$ are the largest $k$ homogeneous blocks of $G$.

Indeed, suppose that $|B_i| \le M$ for some $i \in I(b')$. Note that since the sets $B_i$ are 1-homogeneous blocks, each is contained in some homogeneous block of $G$. Since $|B_i| = 1$ if $i \notin I(b')$, it follows that $t_k \le |B_i| + (|H| - k) \le 2M$. But $G \in \P^{(1)}$, which implies that $t_k \ge N \ge 2M + 1$, so this is a contradiction. Hence $|B_i| \ge M + 1$ for every $i \in I(b')$, as claimed.

Since $|B_i| > 1$ for every $i \in I(b')$, $B_i$ is a homogeneous block for every $i \in I(b')$. Thus, since $|H'| \le k + M$ and $|B_i| \ge M + 1$ for every $i \in I(b')$, $\{B_i : i \in I(b')\}$ are the largest $k$ homogeneous blocks of $G$, which implies that $H' = H_G^{(k)}$, by the definition of $H_G^{(k)}$. We have also shown that $b_G^{(k)}(i) = b'(i) = \infty$ if $i \in I(b')$, and $1$ otherwise, and so $(H',b') = (H_G^{(k)}, b_G^{(k)})$. Since $(H',b')$ were arbitrary, this shows that $(H_G^{(k)}, b_G^{(k)})$ is the unique pair $(H,b)$ such that $G \in \P^*(H,b)$. Call $(H_G^{(k)}, b_G^{(k)})$ the pair `realised by' $G$.

It remains to count how many ordered graphs $G \in \P^{(1)}$ of order $n$ realise a given pair $(H,b) \in \D$. Let $n > k(N + 2M + 1) + M$. Each vertex $i \in V(H) \setminus J(b)$ is assigned one vertex of $G$, and each of the remaining $k$ vertices of $H$ must be assigned at least $N + M + 1 - (|V(H)| - k)$ vertices of $G$, by the definition of $\P^{(1)}$. The remaining $n' = n - k(N + M + 1 - (|V(H)| - k)) - (|V(H)| - k)$ vertices of $G$ may then be assigned arbitrarily to the vertices of $J(b)$. Hence there are exactly ${{n' + k - 1} \choose {k - 1}}$ ordered graphs in $\P^{(1)}_n$ which realise a given pair $(H,b) \in \D$. Therefore
$$|\P^{(1)}_n| = \sum_{(H,b) \in \D} {{n - K(H,b)} \choose {k - 1}},$$ where $K(H,b) = k(N + M + k - 1) - (k - 1)|V(H)| - 1$.

Now, when $k = 1$ we have shown that $|\P_n| = |\D|$ for sufficiently large values of $n$, so the lemma holds in the base case. Let $k \ge 2$, and assume that result holds for all smaller values of $k$. Since $\P^{(4)}$ is a hereditary property of ordered graphs with speed $\Theta(n^{k-2})$, it has speed equal to some polynomial for sufficiently large $n$, by the induction hypothesis. Also, $|\P^{(1)}_n| = \sum_{(H,b) \in \D} {{n - K(H,b)} \choose {k - 1}}$ for sufficiently large $n$, and ${{n - K} \choose {k-1}} = {n \choose {k-1}} - \sum_{i = 1}^K {{n - i} \choose {k - 2}}$ for any $K \in \N$, so  $|\P^{(1)}_n|$ is exactly a polynomial (i.e., it may be written in the form stated in the theorem). Since $\P_n = \P^{(1)}_n \cup \P^{(4)}_n$ and these sets are disjoint, the induction step follows, and hence so does the theorem.
\end{proof}

\section{$\ell$-empty ordered graphs}\label{ellsec}

Let us now return to the general case. We know that $|\P_n| < 2^{n-1}$ for some $n \in \N$, so $|\P_n| < F_{n,\ell+1}$ for some $\ell \in \N$, since $F_{n,\ell+1} = 2^{n-1}$ if $\ell \ge n - 1$. Hence, by Corollary~\ref{summ}, there exist integers $k$ and $\ell$ such that every ordered graph $G \in  \P$ may be partitioned into at most $k$ blocks of consecutive vertices, with each block $\ell$-homogeneous. Before we can deduce the speed of such a property, we need to know more about the ordered graphs induced by these $\ell$-homogeneous blocks. The lemmas in this section will allow us to describe them quite precisely. We start with some definitions, which will make our results much easier to state.

Let $G$ be an ordered graph and $u,v \in V(G)$ (with $u < v$). We say that the pair $u,v$ \emph{separates the edges} of $G$ if for every edge $ij \in E(G)$ with $i < j$, either $j \le u$ or $v \le i$. We will call an ordered graph $G$ \emph{irreducible} if no pair of vertices separate the edges of $G$.

Given any ordered graph $G$, the vertex set of $G$ can be decomposed in a unique way into blocks of consecutive vertices such that each block induces an irreducible subgraph, and there are no edges between different blocks. We call this the \emph{irreducible block decomposition} of $G$, and write $B(G) = (G_1, \ldots, G_m)$ if the irreducible blocks of $G$ induce the ordered graphs $G_1, \ldots, G_m$ in that order.

To be more precise, let $G_1, \ldots, G_m$ be (not necessarily distinct) ordered graphs, let $n_i$ be the order of $G_i$ for each $i \in [m]$, and let $N_i = \sum_{j=1}^i n_j$ for each $i \in [0,m]$. Define $G_1 + \ldots + G_m$ to be the ordered graph with $V(G_1 + \ldots + G_m) = [N_m] = n_1 + \ldots + n_m$, and
\begin{eqnarray*}
E(G_1 + \ldots + G_m) & = & \bigcup_{i = 1}^{m} \; \{ uv : u,v \in [N_{i-1} + 1, N_i]\\
& & \textup{\hspace{1.5cm} and }\;(u-N_{i-1})(v-N_{i-1}) \in E(G_i) \}.
\end{eqnarray*}

Observe that an ordered graph $G$ has the irreducible block decomposition $B(G) = (G_1, \ldots, G_m)$ if and only if $G = G_1 + \ldots + G_m$, and each $G_i$ is irreducible. If $B(G) = (G_1, \ldots, G_m)$ and $|G_i| = n_i$ for each $i \in [m]$, then write $BS(G) = (n_1, \ldots, n_m)$ and call this the \emph{irreducible block sequence} of $G$. Finally, if $G$ is an ordered graph, define $S_n(G)$ to be the number of distinct (i.e., non-isomorphic) induced ordered subgraphs of $G$ of order $n$.

We begin with a simple observation.

\begin{obs}\label{subirred}
Let $n,k \in \N$ with $k \le n$. If $G$ is an irreducible ordered graph of order $n$, then there exists an irreducible ordered subgraph of $G$ of order $k$.
\end{obs}

\begin{proof}
The proof is by induction on $n$. For $n = 1$ and for $k = n$ the result is trivial, so let $n \ge 2$ and $k \le n - 1$. Let $u$ be maximal so that $1u \in E(G)$ and $v$ be maximal so that $2v \in E(G)$. If $u \ge v$ (or $v$ does not exist) then remove vertex 2; if $u < v$ then remove vertex 1. The resulting ordered graph is clearly irreducible, and has order $n-1$, so we are done by induction.
\end{proof}

Our first lemma controls the number of irreducible blocks of size at least $\ell$.

\begin{lemma}\label{KLblocks}
Let $G$ be an ordered graph, and let $k,\ell \in \N$. If there are at least $k$ blocks of order at least $\ell$ in the irreducible block decomposition of $G$, then $S_n(G) \ge F_{n,\ell}$ for each $n \le k$.
\end{lemma}

\begin{proof}
Let $n,k \in \N$ with $n \le k$, and let $G' = G_1 + \ldots + G_k$ be a subgraph of $G$ induced by $k$ of the irreducible blocks of $G$ which are of size at least $\ell$. Thus $G_i$ is irreducible for each $i \in [k]$. For each sequence $(a_1, \ldots, a_t)$ with $t \in \N$, $a_i \in [\ell]$ for each $i \in [t]$ and $\sum_i a_i = n$, choose a subgraph of $G'$ with irreducible block sequence $(a_1, \ldots, a_t)$; such a subgraph exists by Observation~\ref{subirred}, and because $n
\le k$.

These subgraphs are all distinct (since they have distinct irreducible block sequences), so it only remains to count them. It is easy to see that there are exactly $F_{n,\ell}$ sequences $(a_1, \ldots, a_t)$ as described above. Therefore $G$ has at least this many distinct subgraphs of order $n$.
\end{proof}

Using Lemma~\ref{KLblocks}, we can now control the size of the largest irreducible block.

\begin{lemma}\label{bigblock}
Let $k,\ell \in \N$, and let $G$ be an $\ell$-empty ordered graph. If there is a block of size at least $4k\ell$ in the irreducible blocks decomposition of $G$, then $S_n(G) \ge F_{n,\ell}$ for each $n \le k$.
\end{lemma}

\begin{proof}
Let $k,\ell \in \N$ with $n \le k$, and let $G$ be an $\ell$-empty ordered graph with an irreducible block $B$ of size $m \ge 4k\ell$. Let $G'$ be the irreducible subgraph of $G$ induced by $B$, with vertex set $[m]$. We shall find a subgraph $H$ of $G'$ for which $B(H)$ contains $k$ irreducible blocks of size at least $\ell$.

Since $G'$ is irreducible, for each vertex $j \in [m]$ there exists at least one edge $uv \in E(G')$ with $j \in [u,v]$. Thus, for each $i \in [k]$, we may define $$u_i = \min\{u : \exists v, uv \in E(G')\textup{ and }(4i - 3)\ell \in [u,v]\},$$ and $$v_i = \max\{v : \exists u, uv \in E(G')\textup{ and }(4i - 2)\ell \in [u,v]\}.$$

Consider the set of vertices $\{u_1, v_1, \ldots, u_k, v_k\}$. For every $i \in [k]$, $u_i \le (4i-3)\ell$ and $v_i \ge (4i-2)\ell$, so $v_i - u_i \ge \ell$. Also, since $G$ is $\ell$-empty, $u_i \ge (4i-4)\ell + 1$ and $v_i \le (4i-1)\ell - 1$ for every $i \in [k]$. Hence $v_i + \ell < u_{i+1}$ for every $i \in [k-1]$.

Set $A = [u_1,v_1] \cup \ldots \cup [u_k,v_k]$, and let $H = G'[A]$. We claim that $B(H) = (G'[u_1,v_1], \ldots, G'[u_k,v_k])$. Since $v_i - u_i \ge \ell$ for each $i \in [k]$, this implies that $H$ has $k$ irreducible blocks of size at least $\ell$, so, by Lemma~\ref{KLblocks}, it will suffice to prove the lemma.

We must show that $H_i = G'[u_i,v_i]$ is irreducible for each $i \in [k]$, and that there are no edges in $G'$ between $[u_i,v_i]$ and $[u_j,v_j]$ if $i \neq j$. The latter statement follows because $G$ is $\ell$-empty, and $v_i + \ell < u_{i+1}$ for every $i \in [k - 1]$, as observed above. To establish the former statement, suppose that $H_i$ is not irreducible for some $i \in [k]$. Then there must be some consecutive pair $x,y \in [u_i,v_i]$ which separates the edges of $H_i$. By the definitions of $u_i$ and $v_i$, there exist edges $u_iv$ and $uv_i$ in $G'$ with $v \ge (4i-3)\ell$ and $u \le (4i-2)\ell$, so $x,y \in [v,u] \subset [(4i-3)\ell,(4i-2)\ell]$. Since $G'$ is irreducible, there exists an edge $ab \in E(G')$ with $a \le x, y \le b$, and since $x$ and $y$ separate the edges of $H_i$, we must have either $a < u_i$ or $b > v_i$. In either case we have a contradiction, since $u_i$ was chosen to be minimal, and $v_i$ was chosen to be maximal. This contradiction proves that $H_i$ is irreducible, so $B(H) = (G'[u_1,v_1], \ldots, G'[u_k,v_k])$, as claimed.
\end{proof}

The next lemma will allow us to tell which graphs may be induced by arbitrarily many irreducible blocks.

\begin{lemma}\label{ge2}
Let $k \in \N$, let $G_1, \ldots, G_k$ be ordered graphs, and let $G = G_1 + \ldots + G_k$. Suppose that for each $i \in [k]$ there is an integer $m(i)$ such that $G_i$ has at least two irreducible induced subgraphs on $m(i)$ vertices. Then $S_n(G) \ge 2^{n-1}$ for each $n \le k$.
\end{lemma}

\begin{proof}
We prove the lemma by induction on $k$. For $k = 1$ it is trivial, so let $k \ge 2$ and assume it is true for $k - 1$. Let $G = G_1 + \ldots + G_k$ be an ordered graph as described, and let $G' = G_2 + \ldots + G_k$. By the induction hypothesis applied to $G'$, we have $$S_n(G) \ge S_n(G') \ge 2^{n-1}$$ for every $n \le k - 1$. We must show that $S_k(G) \ge 2^{k-1}$.

Let $m = m(1)$, so $G_1$ has two irreducible subgraphs of order $m$, $H_m$ and $H'_m$. By Observation~\ref{subirred}, $G_1$ also has an irreducible subgraph $H_j$ on $j$ vertices for each $j \le m - 1$.

Now, by the induction hypothesis, $G$ has at least $2^{k-j-1}$ induced subgraphs of order $k$ whose left-most irreducible block is $H_j$, for each $1 \le j \le m - 1$. Also $G$ has at least $2^{k-m-1}$ induced subgraphs (of order $k$) whose left-most irreducible block is $H_m$, and at least $2^{k-m-1}$ whose left-most irreducible block is $H'_m$. These induced subgraphs are all distinct (since they have different left-most irreducible blocks). Therefore, $$S_k(G) \ge \sum_{j=1}^{m-1} 2^{k-j-1} + 2 \cdot 2^{k-m-1} = 2^{k-1}.$$
This proves the induction step, and so also the lemma.
\end{proof}

We now define the following five collections of (irreducible) ordered graphs, and two sporadic examples, which will play a central role in our characterisation of properties with speed $p(n)F_{n,k}$. Our reasons for choosing these particular ordered graphs will be made clear by Lemma~\ref{P3P4}. For each $n \in \N$, let
\begin{itemize}
\item $J^{(n)}_1 = K_n$,
\item $J^{(n)}_2$ have vertex set $[n]$ and edge set $E = \{1n\}$ (if $n \ge 2$),
\item $J^{(n)}_3$ have vertex set $[n]$ and edge set $E = \{1i : i \in [2,n]\}$,
\item $J^{(n)}_4$ have vertex set $[n]$ and edge set $E = \{in : i \in [n - 1]\}$,
\item $L^{(n)}$ have vertex set $[n]$ and edge set $E = \{i(i+1) : i \in [n - 1]\}$,
\item $Q_1$ have vertex set $\{1,2,3,4\}$ and edge set $\{13,24\}$,
\item $Q_2$ have vertex set $\{1,2,3,4\}$ and edge set $\{14,23\}$.
\end{itemize}
Also let $\J_k = \{J^{(n)}_i : i \in [4], n \le k\} \cup \{L^{(n)} : n \le k\}$ for $k = 1,2,3$ and let $\J_k = \{J^{(n)}_i : i \in [4], n \le k\} \cup \{L^{(n)} : n \le k\} \cup \{Q_1,Q_2\}$ for each $k \ge 4$. Finally, let $\J = \bigcup_{k \in \N} \J_k$.

\begin{lemma}\label{P3P4}
Let $G$ be a finite irreducible ordered graph, and suppose that $G$ has at most one irreducible ordered subgraph of order $3$, and at most one of order $4$. Then $G \in \J$.
\end{lemma}

\begin{proof}
The result is proved by a simple case analysis, as follows. Let $G$ be an irreducible ordered graph with vertex set $[n]$, and let $t = \max\{|i - j| : ij \in E(G)\}$ be the length of the longest edge in $G$. Suppose that $G$ has at most one irreducible ordered subgraph of order $i$ for $i = 3,4$, and that $G \notin \J$. If $t = 1$ then $G = L^{(n)} \in \J$, since $G$ is irreducible, so $t \ge 2$.

Let $i(i+t) \in E(G)$ be an edge of maximal length in $G$, and suppose first that $i + t < n$. Since $G$ is irreducible, the pair $(i+t,i+t+1)$ does not separate the edges of $G$, so there must be an edge $uv \in E(G)$ with $i < u \le i+t < v$. Now, if $u < i + t$ and $t \ge 3$ then the subgraphs $G[\{i,i+1,i+2,i+t\}]$ and $G[\{i,u,i+t,v\}]$ are distinct (since $iv \notin E(G)$), irreducible subgraphs of $G$, each on 4 vertices, a contradiction. If $u = i + t$ then $G[\{i,i+t,v\}]$ and $G[\{i,i+1,i+t\}]$ are distinct (again since $iv \notin E(G)$), irreducible subgraphs of $G$, each on 3 vertices, another contradiction. So $t = 2$ and $u = i + 1$, i.e., $i(i+2)$ and $(i+1)(i+3) \in E(G)$.

Now, if $j(j+1) \in E(G)$ for some $j \in \{i,i+1,i+2\}$, then $G[i,i+3]$ has at least two irreducible subgraphs on 3 vertices, a contradiction. So if $n = 4$ then $G = Q_1 \in \J$, another contradiction, so $n \ge 5$. Now either $i > 1$ or $i + 3 < n$; assume without loss (by symmetry) that $i+3 < n$. Now, applying to $(i+1)(i+3)$ the same argument that we used for the edge $i(i+t)$, and using the fact that $t = 2$, we obtain $(i+2)(i+4) \in E(G)$. But now $G[\{i,i+1,i+2\}]$ and $G[\{i,i+2,i+4\}]$ are distinct (since $i(i+4) \notin E(G)$) and irreducible subgraphs of $G$, each on 3 vertices, a final contradiction. Therefore $i + t = n$, and similarly it can be proved that $i = 1$.

We have shown that $1n \in E(G)$. Suppose that for some pair $1 < i,j < n$ we have $1i \in E(G)$ and $1j \notin E(G)$. Then $G[\{1,i,n\}]$ and $G[\{1,j,n\}]$ are distinct and irreducible subgraphs of $G$, each on 3 vertices, a contradiction. So $\Gamma(1) = \{n\}$ or $[2,n]$, and similarly $\Gamma(n) = \{1\}$ or $[n - 1]$.

Suppose first that $\Gamma(1) = \{n\}$ and $\Gamma(n) = \{1\}$. If $G[2,n-1] \notin \{E_{n - 2}, K_{n - 2}\}$ then it has an edge $ij$ and a non-edge $uv$, and $G[\{1,i,j,n\}]$ and $G[\{1,u,v,n\}]$ are distinct and irreducible. If $G[2,n-1] = K_{n -2}$ then either $G = Q_2 \in \J$ (if $n = 4$), or $G[\{1,2,n\}]$ and $G[\{2,3,4\}]$ are distinct and irreducible (if $n \ge 5$). If $G[2,n-1] = E_{n-2}$ then $G = J^{(n)}_2 \in \J$. In each case we have a contradiction.

The remaining cases are now easy to deal with. If $\Gamma(1) = [2,n]$ and $\Gamma(n) = \{1\}$ then either $G[2,n-1]$ contains an edge $ij$, in which case $G[\{1,i,j\}]$ and $G[\{1,i,n\}]$ are distinct and irreducible, or $G[2,n-1] = E_{n-2}$, in which case $G = J^{(n)}_3 \in \J$. Similarly if $\Gamma(1) = \{n\}$ and $\Gamma(n) = [n-1]$ then either $G[\{i,j,n\}]$ and $G[\{1,i,n\}]$ are distinct and irreducible or $G = J^{(n)}_4 \in \J$. Finally, if $\Gamma(1) = [2,n]$ and $\Gamma(n) = [n-1]$ then either there is a non-edge $ij$ in $G[2,n-1]$, in which case $G[\{1,i,j\}]$ and $G[\{i,j,n\}]$ are distinct and irreducible, or $G[2,n-1] = K_{n-2}$, in which case $G = J^{(n)}_1 \in \J$. In each case we have a contradiction, so the assumed ordered graph $G$ is impossible, and the proof is complete.
\end{proof}

Combining Lemmas~\ref{ge2} and \ref{P3P4}, we obtain the following result.

\begin{lemma}\label{allJ}
Let $k,m \in \N$, and let $G$ be an ordered graph with irreducible block decomposition $B(G) = (G_1, \ldots, G_m)$. If $|\{i \in [m] : G_i \notin \J\}| \ge k$, then $S_n(G) \ge 2^{n-1}$ for every $n \le k$.
\end{lemma}

\begin{proof}
Let $k,m \in \N$, and let $G$ be an ordered graph whose irreducible block decomposition $B(G) = (G_1, \ldots, G_m)$
satisfies $|\{i \in [m] : G_i \notin \J\}| \ge k$. Let $\{a(1), \ldots, a(k)\} \subset \{i \in [m] : G_i \notin \J\}$, and let $G' = G_{a(1)} + \ldots + G_{a(k)} \le G$.

Now, for each $i \in [k]$ we have $G_{a(i)} \notin \J$, so by Lemma~\ref{P3P4} there is an integer $m(i)$ such that $G_{a(i)}$ has at least two irreducible induced subgraphs on $m(i)$ vertices (and in fact it can be assumed that $m(i) \in \{3,4\}$). Thus, by Lemma~\ref{ge2}, we have $S_n(G) \ge S_n(G') \ge 2^{n-1}$ for every $n \le k$.
\end{proof}

Given an ordered graph $G$ and $k \in \N$, define $G^k = G + \ldots + G$, where $G$ appears $k$ times in the sum.

\begin{lemma}\label{AplusB}
Let $k \in \N$, let $A, B \in \J$ with $A \not\le B$ and $B \not\le A$, and let $G = (A + B)^k$. Then $S_n(G) \ge 2^{n-1}$ for every $n \le k$.
\end{lemma}

\begin{proof}
Note that if $A,B \in \J$, with $A \not\le B$ and $B \not\le A$, then $A + B$ has at least two distinct irreducible ordered subgraphs on $\min(|A|,|B|)$ vertices. Setting $G_i = A + B$ for $1 \le i \le k$, the result follows now immediately by Lemma~\ref{ge2}.
\end{proof}

For each $k,\ell \in \N$, let $\J(k,\ell)$ be the following collection of ordered graphs. Let $G \in \J(k,\ell)$ if and only if there exist an (ordered) collection of $s \le k$ ordered graphs $(A_1, \ldots, A_s)$ satisfying the following conditions.
\begin{itemize}
\item $G = A_1 + \ldots + A_s$, and
\item For each $i \in [s]$, there exists an (ordered) collection of ordered graphs $(B^{(i)}_1, \ldots, B^{(i)}_{t(i)})$, satisfying
\begin{itemize}
\item $A_i = B^{(i)}_1 + \ldots + B^{(i)}_{t(i)}$,
\item $B^{(i)}_j \in \J_\ell$ for each $j \in [t(i)]$, and
\item $B^{(i)}_j \le B^{(i)}_{j'}$ or $B^{(i)}_{j'} \le B^{(i)}_j$ for each pair $j,j' \in [t(i)]$.
\end{itemize}
\end{itemize}
We call a collection $(A_1, \ldots, A_s)$ satisfying these conditions a \emph{$(k,\ell)$-witness set} for $G$.

\begin{lemma}\label{Jkl}
Let $k,\ell,m \in \N$ with $k \ge 50m\ell^2 + 1$, let $G$ be an ordered graph, and suppose that $G \in \J(k,\ell)$ but $G \notin \J(k-1,\ell)$. Then $S_n(G) \ge 2^{n-1}$ for every $n \le m$.
\end{lemma}

\begin{proof}
Let $k,\ell,m \in \N$, with $k \ge 50m\ell^2 + 1$, let $G$ be an ordered graph, and suppose that $G \in \J(k,\ell)$ but $G \notin \J(k-1,\ell)$. Let $\{A_1, \ldots, A_k\}$ be a $(k,\ell)$-witness set for $G$. Since $k$ was chosen to be minimal, for each $i \in [k-1]$ there must exist irreducible blocks $D_i \subset A_i$ and $D_i' \subset A_{i+1}$ such that $G[D_i] \not\le G[D'_i]$ and $G[D'_i] \not\le G[D_i]$; otherwise for some $i$ the collection $(A_1, \ldots, A_{i-1}, A_i + A_{i+1}, A_{i+2}, \ldots, A_k)$ would be a $(k-1,\ell)$-witness set for $G$.

Consider the multiset of pairs $\D = \{(D_{2i-1},D_{2i-1}') : 2i \le k\}$. Since $k \ge 50m\ell^2 + 1$, $|\D| \ge 25m\ell^2$. Now, note that there are fewer than $5\ell$ ordered graphs in $\J_\ell$, and that $D_{2i-1}$ and $D_{2i-1}' \in \J_\ell$ for each $i$ with $2i - 1 \le k$. Thus, by the pigeonhole principle, there exist at least $m$ copies of some pair $(D,D')$ in $\D$. Therefore, $H = (D + D')^m \le G$, so by Lemma~\ref{AplusB}, $S_n(G) \ge S_n(H) \ge 2^{n-1}$ for every $n \le m$.
\end{proof}

Finally, we make the following observation.

\begin{obs}\label{symmdiff}
Let $k,\ell,n \in \N$, let $G$ be an ordered graph on $[n]$, and suppose that $G$ can be partitioned into $k$ blocks of consecutive vertices, with each block $\ell$-homogeneous. Then there exists an ordered graph $H$ on $[n]$ with at most $k$ homogeneous blocks, such that $G \triangle H$ is $\ell$-empty.
\end{obs}

\begin{proof}
Let the $\ell$-homogeneous blocks of $G$ be $B_1, \ldots, B_k$, and let $x \in B_i$ and $y \in B_j$ with $i,j \in [k]$. Either all or none of the edges of length at least $\ell$ between $B_i$ and $B_j$ are in $G$. Let $xy \in E(H)$ if and only if all of these edges are in $G$. (So if there are no edges of length at least $\ell$, then $xy \notin E(H)$.) Then $G \triangle H$ is $\ell$-empty.
\end{proof}

\section{The structure of a property with speed $p(n)F_{n,\ell}$}\label{Fnksec}

We can now deduce the structure of every ordered graph $G \in \P$, if $\P$ is a hereditary property whose speed satisfies $|\P_n| < F_{n,\ell+1}$ for some $n,\ell \in \N$.

\begin{thm}\label{fibstruc}
Let $n,\ell \in \N$, let $\P$ be a hereditary property of ordered graphs, and suppose that $|\P_n| < F_{n,\ell+1}$. Then there exist $k,k' \in \N$ such that every ordered graph $G \in \P$ is of the form $G = H \triangle J$, where $H$ is an ordered graph with at most $k$ homogeneous blocks, and $J \in \J(k',\ell)$.
\end{thm}

\begin{proof}
Let $n,\ell \in \N$, let $\P$ be a hereditary property of ordered graphs, and suppose that $|\P_n| < F_{n,\ell+1} \le 2^{n-1}$. By Corollary~\ref{summ}, there exists an integer $K \in \N$ such that every ordered graph $G \in \P$ may be partitioned into at most $K$ blocks of consecutive vertices with each block $\ell$-homogeneous.

Let $G \in \P$, let $B$ be an $\ell$-homogeneous block of $G$, and let $F = G[B]$. Suppose that $F$ is $\ell$-empty. By Lemma~\ref{KLblocks}, $F$ has at most $n-1$ irreducible blocks of size at least $\ell + 1$, and by Lemma~\ref{bigblock}, $F$ has no irreducible block of size at least $4n(\ell + 1)$, since $S_n(F) \le |\P_n| < F_{n,\ell+1}$. Therefore, by deleting a set of edges which span at most $4n^2(\ell + 1)$ vertices of $F$, we can obtain an ordered graph $F'$ in which each irreducible block has size at most $\ell$.

Now, by Lemma~\ref{allJ}, at most $n - 1$ of these irreducible blocks are not in $\J$, since $|\P_n| < 2^{n-1}$. So by deleting a set of edges which span at most $(n-1)\ell$ vertices from $F'$, we can obtain an ordered graph $F''$ in which each irreducible block is in $\J_\ell$.

Finally, let $s \in \N$ be minimal such that $F'' \in \J(s,\ell)$. By Lemma~\ref{Jkl} we have $s \le 50n\ell^2$, since $|\P_n| < 2^{n-1}$.

We have shown that if $B$ is an $\ell$-homogeneous block of $G$ which induces an $\ell$-empty graph $F$, then by deleting a set of edges which span at most $4n^2(\ell + 1) + (n-1)\ell < 9n^2\ell$ vertices of $F$, we can obtain an ordered graph in $\J(50n\ell^2,\ell)$. By symmetry, if $B$ induces an $\ell$-complete graph $F^*$, then by adding a set of non-edges which span at most $9n^2\ell$ vertices of $F^*$, we can obtain a graph whose complement is in $\J(50n\ell^2,\ell)$. For each $\ell$-homogeneous block $B_i$ of $G$, choose such a collection of edges, $E_i$. Let $H'$ be the ordered graph on $[N] = V(G)$, with edge set $\bigcup_i E_i$, and note that $H'$ is $\ell$-empty.

Now, every ordered graph in $\P$ may be partitioned into $K$ or fewer $\ell$-homogeneous blocks, $B_1 < \ldots < B_K$. Thus, by Observation~\ref{symmdiff}, there exists an ordered graph $H''$ on $[n]$ with at most $K$ homogeneous blocks, such that $G \triangle H''$ is $\ell$-empty.

Finally, we need to remove edges of $G \triangle H''$ between different $\ell$-homogeneous blocks of $G$, so let $H'''$ have vertex set $[N]$ and edge set $E(H') \cup \{uv \in E(G \triangle H'') : u \in B_i, v \in B_j, i \neq j\}$.

Let $H = H'' \triangle H'''$. We claim that $H$ has at most $18K n^2 \ell + 4K\ell + K$ homogeneous blocks, and that $G \triangle H \in \J(k',\ell)$ if $k' \ge 50Kn\ell^2$. The first statement follows because fewer than $9Kn^2\ell + 2K\ell$ vertices of $H'''$ have non-zero degree, and $H''$ has at most $K$ homogeneous blocks. To prove the second statement, note that $G \triangle H = (G \triangle H'') \triangle H'''$, so $G \triangle H = A_1 + \ldots + A_K$, where $A_i \in \J(50n\ell^2,\ell)$ for each $i \in [K]$. It is now easy to see that $G \triangle H \in \J(k',\ell)$ if $k' \ge 50Kn\ell^2$.

Thus, letting $k = 18Kn^2\ell + K + 4K\ell$ and $k' = 50Kn\ell^2$, we have $G = H \triangle J$ for some $J \in \J(k',\ell)$.
\end{proof}

\section{Proof of Theorem~\ref{order}}\label{pfsec}

We shall use the following easy fact about Fibonacci numbers.

\begin{obs}\label{fibadd}
Let $\ell, m, n \in \N$. Then $F_{m + n, \ell} \ge F_{m,\ell} \cdot F_{n,\ell}$.
\end{obs}

\begin{proof}
Let $\ell \in \N$. We use induction on $m + n$. We have $F_{2,\ell} = 2 > 1 = F_{1,\ell} \cdot F_{1,\ell}$, so the result holds for $m + n = 2$. Now, let $m + n > 2$, and assume the result is true for all smaller values of $m + n$. Then
\begin{align*}
F_{m + n, \ell} & = \; F_{m+n-1,\ell} + \ldots + F_{m+n-\ell,\ell} \\
& \ge \; F_{m,\ell} \left( F_{n-1,\ell} + \ldots + F_{n-\ell,\ell} \right) \; = \; F_{m,\ell} \cdot F_{n,\ell}.
\end{align*}
So the induction step holds, and the observation is proved.
\end{proof}

\noindent The following lemmas provide the final piece of the jigsaw.

\begin{lemma}\label{sizeJ1}
Let $\ell, n \in \N$. Then $|\J(1,\ell)_n| = O(F_{n,\ell})$ as $n \to \infty$.
\end{lemma}

\begin{proof}
Let $\ell, n \in \N$, and let $G \in \J(1,\ell)_n$. Then for some $t \in \N$ we have $G = B_1 + \ldots + B_t$, with $B_j \in \J_\ell$ for each $j \in [t]$, and $B_i \le B_j$ or $B_j \le B_i$ for each pair $i,j \in [t]$. This is the irreducible block decomposition $B(G)$ of $G$, and so is clearly unique. Note also that there is a unique `largest' ordered graph $B$ in $B(G)$, i.e., $B_j \le B$ for each $j \in [t]$, and $B_i = B$ for some $i \in [t]$.

How many such ordered graphs $G$ are there? Partition $\J(1,\ell)$ as follows: for each $B \in \J_\ell$, let $\J(B) = \{G \in \J(1,\ell) : B$ is the largest irreducible graph in $B(G)\}$. Each ordered graph $B \in \J_\ell$ has order at most $\ell$, so there are only a bounded number of them. Hence we will be done if we can prove that $|\J(B)_n| = O(F_{n,\ell})$ for every $B \in \J_\ell$.

But this is now easy, since every ordered graph in $\J(B)$ is a subgraph of $B^m$ for some sufficiently large $m$. Now simply observe that for each $B \in \J_\ell$, $B$ has exactly one irreducible ordered subgraph of order $n$ for each $n \le |B|$, and it follows by a simple induction on $n$ that $|\J(B)_n| = F_{n,|B|} \le F_{n,\ell}$ for every $n \in \N$.
\end{proof}

Using Lemma~\ref{sizeJ1}, we can now give an upper bound on $|\J(k,\ell)_n|$ for all $k$, $\ell$ and $n \in \N$.

\begin{lemma}\label{sizeJk}
Let $k, \ell, n \in \N$. Then $|\J(k,\ell)_n| = O(n^{k-1}F_{n,\ell})$ as $n \to \infty$.
\end{lemma}

\begin{proof}
Let $k, \ell, n \in \N$. If $G \in \J(k,\ell)_n$, then $G = A_1 + \ldots + A_k$, with each $A_i \in \J(1,\ell) \cup \emptyset$ (where $\emptyset$ here denotes the ordered graph with $|G| = 0$). Let $a(i) = |A_i|$ for each $i \in [k]$, so $a_i \in \N \cup \{0\}$.

It follows from this that an ordered graph $G \in \J(k,\ell)_n$ is determined by a sequence $(a(1), \ldots, a(k))$, with $a(i) \in \N \cup \{0\}$ for each $i \in [k]$, and $\sum_i a_i = n$; and a sequence of ordered graphs $(A_1, \ldots, A_k)$, with $A_i \in \J(1,\ell)_{a(i)}$ for each $i \in [k]$. There are $O(n^{k-1})$ such sequences of integers, and by Observation~\ref{fibadd} and Lemma~\ref{sizeJ1} there are $$\prod_{i=1}^k O(F_{a(i),\ell}) = O(F_{a(1) + \ldots + a(k), \ell}) = O(F_{n,\ell})$$ such sequences of ordered graphs. The result follows immediately.
\end{proof}

\begin{proof}[Proof of Theorem~\ref{order}]
Let $\P$ be a hereditary property of ordered graphs, and suppose that $|\P_n| < 2^{n-1}$ for some $n \in \N$. In particular, let $m \in \N$ satisfy $|\P_m| < 2^{m - 1}$. Since $F_{n,\ell} = 2^{n-1}$ for every $n \le \ell$, it follows that $|\P_m| < F_{m,m}$. Let $\ell \in \N$ be the minimal integer such that $|\P_n| < F_{n,\ell+1}$ for some $n \in \N$; we have shown that $\ell \le m - 1$, so such an integer $\ell$ exists. Note that by the definition of $\ell$, $|\P_n| \ge F_{n,\ell}$ for every $n \in \N$.

Suppose first that $\ell = 1$, so $|\P_n| < F_n$ for some $n \in \N$. Theorem~\ref{poly} and Lemma~\ref{geton} then imply that either case (a) or case (b) of the theorem holds. So let $\ell \ge 2$, and apply Theorem~\ref{fibstruc} to $\P$. By the theorem, there exist integers $k,k' \in \N$ such that every ordered graph $G \in \P$ may be written as $G = H \triangle J$, where $H$ has at most $k + 1$ homogeneous blocks, and $J \in \J(k'+1,\ell)$. By Lemma~\ref{blocks} there are $O(n^k)$ such ordered graphs $H$ on $n$ vertices, and by Lemma~\ref{sizeJk} there are $O(n^{k'}F_{n,\ell})$ such ordered graphs $J$ on $n$ vertices. Hence $|\P_n| = O(n^{k + k'}F_{n,\ell})$.
\end{proof}

\section{Further problems}\label{probsec}

Theorems~\ref{order} and \ref{highord} restrict the possible speeds of a hereditary property of ordered graphs if the speed is at most $2^{n-1}$, or at least $2^{cn^2}$ for some $c > 0$. There are many obvious questions remaining in the large gap between these ranges, and in this section we shall discuss some of these.

For hereditary properties of both labelled graphs and permutations, there is a jump from exponential speed (speed $c^n$ for some constant $c > 0$) to factorial speed (speed $n^{cn}$ for some constant $c > 0$) (see \cite{BBW1} and \cite{MT}). As we have seen, ordered graphs generalize both of these types of structure, so it is natural to ask whether a similar jump occurs for ordered graphs. In \cite{BBMpar} it was proved that such a jump does occur for hereditary properties of ordered graphs in which every component is a clique, for monotone properties of ordered graphs (properties closed under taking arbitrary (i.e., not necessarily induced) ordered subgraphs), and for hereditary properties of ordered graphs not containing arbitrarily large complete, or complete bipartite ordered graphs. (The first two of these results have been proved independently by Klazar and Marcus~\cite{KM}.) It was also conjectured that the same jump holds for arbitrary hereditary properties of ordered graphs.

Even assuming a positive answer to the conjecture, one is still left with the problem of determining the possible exponential speeds. In particular, we have the following questions.

\begin{conj}\label{ctothen}
If $\P$ is a hereditary property of ordered graphs, and $|\P_n| < c^n$ for some $c \in \RR$ and every $n \in \N$, then $\ds\lim_{n \to \infty} (|\P_n|)^{1/n}$ exists.
\end{conj}

\begin{prob}\label{bases}
Let $\S = \{c \in \RR :$ there is a hereditary property of ordered graphs $\P$ with $\ds\lim_{n \to \infty} (|\P_n|)^{1/n} = c\}$. Determine the set $\S$.
\end{prob}

Theorem~\ref{order} solves Problem~\ref{bases} in the case $\ds\liminf_{n \to \infty} (|\P_n|)^{1/n} < 2$. The following corollary is immediate from the theorem.

\begin{cor}
Let $\S$ be as defined in Problem~\ref{bases}, and let $A = \{x : x$ is the largest real root of the polynomial $x^{k+1} = x^k + x^{k-1} + \ldots + 1$ for some $k \in \N\}.$ Then $\ds\S \cap [0,2] = \{0,2\} \cup A.$
\end{cor}

Arratia~\cite{Arratia} proved Conjecture~\ref{ctothen} for principal hereditary properties of permutations. The following easy result (which uses basically the same method) proves another special case of the conjecture.

\begin{thm}\label{fekete}
Let $G_1, G_2, \ldots$ be a sequence of ordered graphs, and suppose that either every $G_i$ is irreducible, or every $\overline{G_i}$ is irreducible. Let $\P = \{G : G$ is an ordered graph, and $G_i \not\le G$ for every $i \in \N\}$. Then either $\ds \lim_{n \to \infty} (|\P_n|)^{1/n}$ exists, or $\ds \liminf_{n \to \infty} (|\P_n|)^{1/n} = \infty$.
\end{thm}

\begin{proof}
We claim that for every pair of integers $m,n$, $$|\P_{n+m}| \ge |\P_n| \cdot |\P_m|.$$ Assume that every $G_i$ is irreducible (the other case can be dealt with similarly). Let $F_1 \in \P^n$ and $F_2 \in \P^m$, and let $i \in \N$. By the definition of $\P$ we have $G_i \not\le F_1$ and $G_i \not\le F_2$, and so, since $G_i$ is irreducible, $G_i \not\le F_1 + F_2$. This holds for every $i \in \N$, so $F_1 + F_2 \in \P^{n+m}$. This proves the claim.

Now, Fekete's Lemma~\cite{Fek} states that if $a_1, a_2, \ldots \in \RR$ satisfy $a_m + a_n \ge a_{m+n}$ for all $m,n \ge 1$, then $\ds\lim_{n \rightarrow \infty} \displaystyle\frac{a_n}{n}$ exists and is in $[-\infty,\infty)$. Applying this lemma to the sequence $-\log(|\P_n|)$ gives the result.
\end{proof}

There has been a large volume of work done on the possible exponential speeds of principal hereditary properties of permutations (see for example \cite{Bonabook}). Until recently all such known speeds were of the form $k^{(1 + o(1))n}$, with $k \in \N$, but a non-integer base was found by B\'ona~\cite{Bona2}, who proved that the property of all permutations avoiding $12453$ has speed $(9 + 4\sqrt{2})^{(1 + o(1))n}$. However, there are no known hereditary properties with speed $c^{(1 + o(1))n}$ and $c$ transcendental. We have been unable to find even a hereditary property of ordered graphs with transcendental base, but the following simple construction shows that for properties of ordered graphs, irrational bases are much easier to come by than in the more restrictive principal permutation property setting.

\begin{thm}\label{somebases}
Let $k \in \N$, and let $a(0) \le \ldots \le a(k)$ with $a(i) \in \N$ for each $i$. Let $c$ be the largest real root of the polynomial $x^{k+1} = \sum_{i = 0}^k a(i)x^i$. Then there exists a hereditary property of ordered graphs $\P$ with $|\P_n| = c^{(1+o(1))n}$.
\end{thm}

\begin{proof}
Let $k \in \N$, and let $a(0) \le \ldots \le a(k)$ with $a(i) \in \N$ for each $i \in [0,k]$. We shall define a particular infinite ordered graph $G$, and let $\P$ be the property of ordered graphs consisting of all (finite, order-preserving) subgraphs of $G$.

For each $\ell \in [k]$ let $K_\ell$ denote a copy of the complete ordered graph on $\ell$ vertices, and let $H = K_{k+1}^{a(0)} + K_k^{a(1)-a(0)} + \ldots + K_1^{a(k)-a(k-1)}$. Let $BS(H) = (b(1), \ldots, b(a(k)))$, for each $j \in [a(k)]$ let $d(j) = \sum_{i=1}^{j-1} b(i)$, and let $d = \sum_{i=1}^{a(k)} b(i) = |H|$. Let $G$ be an infinite ordered graph with vertex set $\N$, satisfying the following conditions:
\begin{enumerate}
\item[$(i)$] $G - [a(k)] = H + H + \ldots$, and
\item[$(ii)$] for each $i \in [a(k)]$ and each $a(k) < j$, $ij \in E(G)$ if and only if $j - a(k) \in [d(i)+1,d] \pmod d$.
\end{enumerate}

Let $\P$ be the collection of all (finite, order-preserving) induced subgraphs of $G$. It is easy to see that $\P$ is a hereditary property of ordered graphs. Let $T_n$ be the sequence of integers defined by $T_n = 0$ if $n < 0$, $T_0 = 1$, and $T_{n+1} = \sum_{t=0}^k a(k-t)T_{n-t}$ for every $n \ge 0$. We claim that $|\P_n| = \Theta(T_n)$.

Let us first show that $|\P_{n+a(k)}| \ge T_n$ for every $n \in \N$. Indeed, let the first (leftmost) $a(k)$ vertices of $G$ be denoted $A$, and let's consider only those subgraphs of $G$, on $n + a(k)$ vertices, which include all of $A$. Such an ordered graph consists of $A$, and then a sequence of cliques, each connected to a subset (in fact an initial segment) of $A$. If the clique has size $t$, then there are exactly $a(k+1-t)$ choices for this subset. It now follows easily by induction on $n$ that there are exactly $T_n$ such ordered subgraphs.

Now, note that any graph in $\P_n$ may be obtained by first taking an ordered subgraph $J$ of $G$ on $n + i$ vertices (with $0 \le i \le a(k)$) containing all of $A$, and then removing $i$ vertices of $A$. Since $T_n$ is increasing, there are at most $T_n$ choices for $J$ (for a given $i$), and hence $|\P_n| \le 2^{a(k)}T_n$. This proves that $|\P_n| = \Theta(T_n)$, and the theorem follows.
\end{proof}

An \emph{accumulation point from below} of a set $S \subset \RR$ is a point $c \in \RR$ such that for every $\eps > 0$, $S \cap (c - \eps, c) \neq \emptyset$. Let $A^1(S)$ denote the accumulation points from below of $S$, and for each $n \in \N$, let $A^{n+1}(S)$ denote the accumulation points from below of the set $A^n(S)$. We call the set $A^n(S)$ the \emph{degree $n$ accumulation points from below}. Using Theorem~\ref{somebases}, we can obtain the following result about the accumulation points of $\S$.

\begin{cor}\label{accum}
For each $2 \le n \in \N$ we have $n \in A^{n-1}(\S)$. In other words, $n$ is a degree $n-1$ accumulation point from below of $\S$.
\end{cor}

\begin{proof}
We claim that the result holds even if we consider only the family of properties described in Theorem~\ref{somebases}, and prove this claim by induction on $n$. For $n = 2$, consider the sequences $(1)$, $(1,1)$, $(1,1,1)$, and so on, and apply Theorem~\ref{somebases}. This gives a sequence of constants $c_i \in \S$ with $c_i \to 2^-$ as $i \to \infty$ (note that in fact $c_i = \ds\lim_{n \to \infty} (F_{n,i})^{1/n}$), so the claim is true for $n = 2$.

So let $3 \le n \in \N$ and assume the claim holds for $n - 1$. Consider the collection $\A$ of sequences which proved the result for $n-1$. Now add 1 to each entry of each sequence in $\A$, to obtain the collection $\A'$. The sequences in $\A'$ all still satisfy the conditions of Theorem~\ref{somebases}, so we may apply the theorem to them. It follows that $n$ is a degree $n-2$ accumulation point from below of $\S$.

Now, let $\omega = (a(1), \ldots, a(k)) \in \A'$, and suppose applying Theorem~\ref{somebases} to $\omega$ shows that $c \in \S$. We sub-claim that $c \in A^1(\S)$; since $\omega$ was arbitrary, this will suffice to prove the claim. Note that $a(k) \ge 2$, consider the sequences $(a(1), \ldots, a(k), 1)$, $(a(1), \ldots, a(k), 1, 1)$, and so on, and apply Theorem~\ref{somebases}. The theorem gives a sequence of constants $c_i \in \S$ with $c_i \to c^-$ as $i \to \infty$, so $c \in A^1(\S)$ as sub-claimed. Thus $n \in A^{n-2}(A^1(\S)) = A^{n-1}(\S)$, and the induction step is complete. The result follows immediately.
\end{proof}

Of course it is not necessary to use copies of the complete graph in the proof of Theorem~\ref{somebases} -- one could use the irreducible graphs from any finite hereditary property of ordered graphs. However it does not appear that arbitrary sequences are possible (at least using this method). The following conjecture is motivated by Theorem~\ref{somebases} and by Lemma~\ref{P3P4}.

\begin{conj}\label{jumpat2}
The smallest $c \in \S$ with $c > 2$ is the largest real root of the polynomial $x^5 = x^4 + x^3 + x^2 + 2x + 1$, and is approximately $2.03$.
\end{conj}

An example of a hereditary property with this speed is the following. Let $\P$ consist of all ordered graphs in which each irreducible block is either $J_2^{(n)}$ with $n \le 5$, or $Q_1$. It is easy to show that $\P$ is hereditary and has the desired speed. The following conjecture is a much more general version of Conjecture~\ref{jumpat2}. It says that there is a jump everywhere!

\begin{conj}
For every $c \in \RR$, there exists an $\eps = \eps(c)$ such that $\S \cap (c,c+\eps) = \emptyset$. In particular, $\S$ has no accumulation points from above.
\end{conj}

All the members of $\S$ we have found are either integers or algebraic irrationals. Our final conjecture says that all members of $\S$ are of one of these two types.

\begin{conj}
Every $c \in \S$ is either an integer or an irrational algebraic number, so $\S \subset \Z \cup (\AA \setminus \QQ)$.
\end{conj}

Finally, a bipartite ordered graph is a bipartite graph together with a linear order on each part. Note that bipartite ordered graphs are equivalent to $\{0,1\}$-matrices, and that a bipartite ordered graph may be mapped to an ordered graph by placing one part to the left of the other. It is easy to see that for bipartite ordered graphs, structures of Type 3 do not occur. We thus obtain the following corollaries to the proof of Theorem~\ref{order}.

\begin{cor}\label{bipkey}
There exists a function $g: \N \to \N$ such that if $k \in \N$ and $G$ is a bipartite ordered graph containing no $k$-structure of Type 1 or Type 2, then $G$ can be partitioned into at most $g(k)$ homogeneous blocks.
\end{cor}

\begin{cor}
If $\P$ is a hereditary property of bipartite ordered graphs, then either
\begin{enumerate}
\item[(a)] $|\P_n| = \ds\sum_{i = 0}^k a_i {n \choose i}$ for some $k \in \N$, $a_0, \ldots, a_k \in \Z$ and all sufficiently large $n$, or\\[-1ex]
\item[(b)] $|\P_n| \ge 2^{n-1}$ for every $n \in \N$.\\[-1ex]
\end{enumerate}
\end{cor}

\end{document}